Student Sense-Making in Post-Secondary Introductory Proof Courses:

An Argument for and Outline of a Methodological Approach

Bolanle Salaam
Department of Mathematics and Science Education
University of Georgia

Submitted January 11, 2022

Contents

# Chapter 1: Introduction

**Problem Statement**

The act of proof construction can be viewed as a problem-solving task (Weber, 2001). Weber (2001) defines the nature of problem-solving tasks in the following way:

> In problem-solving tasks, the problem solver typically is presented with an initial state and is asked to perform a sequence of actions that will transform the initial state into a desired goal state. In constructing the proof of a statement, the prover is given an initial set of assumptions and is asked to derive a sequence of inferences (e.g. recall definitions and apply theorems) which conclude with the statement to be proven. (pp. 110-111)

The way the participant faced with the problem-solving task chooses to interact with the task can be described in terms of the theory of commognition[1] (Sfard, 2008, 2009, 2012, 2020). Upon receiving[2] the problem-solving task, the participant may notice a set of *signifiers* (Sfard, 2008), defined as words, symbols, gestures, tactile, or artifacts that function as nouns in the utterances of discourse participants. The participant may then associate *realizations*, defined as a procedure that pairs a *signifier* with another object (Sfard, 2008), to the signifiers the participant notices within the task, leading them to enact a particular *routine*.

Due to the centrality of proof construction to the advanced study of mathematics, proof is a major focus of research on undergraduate mathematics. However, proof construction persists as a topic of difficulty for high school students in mathematics, pre-service mathematics teachers, and mathematics majors (Edwards & Ward, 2004; Moore, 1994; Schoenfeld, 1991; Weber, 2001). The discursive competence, though elsewhere regarded as content knowledge, thought to be required to be "successful" in proof construction is typically split into two overlapping categories: discourse about *proof content*, where the phrase *proof content* is often used to refer to the common methods and techniques involved in the act of proving, and discourse about non-proof mathematical content. The individual discourses that collectively define the curriculum, often referred to elsewhere as *concepts*, in transition to proof courses are widely treated as entities in their own right in proof literature. When viewed in the discursive context, the term *concept* can be operationalized by defining it as a word or a signifier together with its discursive use (Sfard, 2008). Therefore, in this view, a term or signifier that is labeled as a *concept* (e.g. proof, function, ring, field, addition, isomorphism, etc.) does not exist independently of the community or individual that determines the signifier or term's discursive uses. Succinctly, teaching a *concept* is a matter of describing, coordinating, and exemplifying the acceptable uses of a term as defined by the relevant community of discourse.

---

[1] The terminology used to develop the necessary background for this study is extensive. A glossary of commognitive terms is linked for the benefit of the reader. See: Thurston, A., & Wall, K. (Eds.). (2012). Developing mathematical discourse-Some insights from communicational research. [Special Issue]. *International Journal of Educational Research*, *51*, 155-158. doi: 10.1016/j.ijer.2012.01.003

[2] To emphasize that the community of people who can engage in mathematical problem-solving is inclusive of the blind, non-verbal, Deaf-blind, and wider Deaf communities in addition to all populations that do not rely on seeing, writing, speaking, and/or hearing for communication, I use the term *receiving* to indicate communication of a task. The population of focus in this study consists of the simultaneous members of the hearing, writing, seeing, verbal, and English-speaking communities.

The non-proof concept of function plays a critical role not only in student proof construction, but to the mathematical development of mathematics majors and the wider student body. Function, which appears in most mathematics courses and is emphasized in middle and high school algebra, is sometimes referred to as a *threshold concept*. Threshold concepts can be described as a bridge to a view of subject matter that was previously inaccessible to the learner. (Meyer & Land, 2005). Once a student crosses this bridge, they have radically changed "understandings" of the concept and related subject area that are irreversible (Meyer & Land, 2005; Pettersson, 2012). Crossing this threshold can be complicated by the difficulty students may face with prior relevant concepts, often blocking students from further advancement in the subject matter. In terms of discourse, it can be said that failing to *individualize*, defined as the process that enables a person to individually act when they previously could only act in concert with others (Sfard, 2008), the threshold concept into the participant's current discourse can prevent them from accessing more advanced levels of the discourse. Function is notably emblematic of this situation (Pettersson, 2012). Thus, it is unsurprising that special types of functions would also be difficult for students to individualize without the appropriate introduction to foundational aspects of function discourse. For example, Breidenbach, Dubinsky, Hawks & Nichols argue that the individualization of injections (one-to-one) and surjections (onto) is likely bolstered by a student's individualization of the threshold concept of functions as a "dynamic transformation of objects according to some repeatable means that, given the same original object, will always produce the same transformed object" (1992, p. 251). Injections and and surjections will be the focal non-proof content in this study.

"Understanding" of injections and surjections seemingly requires students that receive hearing, writing, seeing, and verbally focused instruction to engage with injection and surjection *word use* (e.g. definitions), read and use symbolic *visual mediators* associated with these word uses, and eventually, determine when *applicability conditions* have been met for enacting a particular routine. Within mathematics education literature, a specific instance of applicability conditions and the routine that is triggered when a student encounters a task situation is most similar to what Tall & Vinner (1981) termed *evoked concept image* on the part of the student, whereas the *metarules* an individual appeals to when working with these words along with all appropriate instances of word use, applicability conditions leading to routines, and associated visual mediators would most closely resemble what Tall & Vinner call *concept definition* and *concept image,* defined as "the total cognitive structure that is associated with the concept, which includes all the mental pictures and associated properties and processes" (Tall, 1992; Tall & Vinner, 1981). However, the typical manner of hearing, writing, seeing, and verbally focused instruction regarding one-to-one and onto functions (and consequently, the idea of injective and surjective maps) has features that can be exploited to construct proofs about injective and surjective functions, even in the absence of a student conforming to the *norms* of the wider mathematical community. In the presence of a seemingly "correct" solution to the task, a lack of agreement between student metarules and the norms common to the mathematical community can therefore go undetected by a course instructor.

The potential for students to be able to construct proofs in response to tasks that involve using the definitions of injectivity and surjectivity even when their reproducible discourse on or *concept image* of function may be unformed, or loosely formed, makes injections and surjections ideal for studying student sense-making when confronted with proof tasks. Furthermore, injections are of particular interest as they are introduced to students in the United States by the latter part of their K-12 experience, giving students experience with function and the property of being one-to-one well before their matriculation in university. This allows for an examination of changes in student *thinking*, operationalized as the individualization of interpersonal communication (Sfard, 2008), where it can be reasonably asserted that students will not be "blank slates" in their function discourse upon entering the introductory proof classroom, as this is a rarity when teaching any topic. The different starting points of students is therefore an important feature of, rather than an analytical inconvenience to, the study of changes in students' function discourse over time (Vinner & Dreyfus, 1989, p. 365).

The problem students within the study population of focus (referred to heretofore as undergraduate students belonging to the hearing, writing, seeing, and verbal communities simultaneously) face when encountering injective and surjective functions in introductory proof courses is multi-faceted. Students may: (1) struggle with the notion of function in general, and this struggle may carry over from prior educational settings; (2) have trouble extending prior "understandings" of function to higher education settings due to differences in function representation; (3) struggle with what it "means" to prove a mathematical statement; (4) struggle with proof methods commonly taught in introduction to proof courses; or (5) struggle with deciding when and how to execute known proof strategies and methods. While student difficulty with proof (e.g., Dreyfus, 1999; Edwards & Ward, 2004; Harel & Sowder, 1998; Moore, 1994; Weber, 2001) and student difficulty with the notion of function (e.g., Breidenbach, Dubinsky, Hawks & Nichols, 1992; Dubinsky & Wilson, 2013; O'Shea, Breen & Jaworski, 2016; Selden & Selden, 1992; Tall & Vinner, 1981; Vinner, 1983; Vinner & Dreyfus, 1989) are well documented separately in the literature, Thoma & Nardi (2020) emphasize that research on student engagement at the collegiate level with injective and surjective functions is scarce (e.g., Nardi, Ryve, Stadler, & Viirman, 2014). While some initial advances have been made in this area, rarer still is the collegiate-level study of injective of surjective functions within the specific context of proof (e.g., Thoma & Nardi, 2018, 2020; Thoma, (2018); Bansilal, Brijlall, & Trigueros, 2017).

The reasons for this scarcity are unclear, but I suspect that: (1) injection and surjection proofs are not seen as an area of substantial difficulty for students by researchers and instructors of introductory proof courses or (2) that examining injections and surjections in proof courses is a daunting task due to profound student difficulty with both proof and function individually. Even in the presence of a variety of definitions and beliefs about proof and what constitutes proof within the mathematics education literature (Healy & Hoyles, 2000; Harel & Sowder, 1998), proof can be thought of as a mathematical concept itself across these definitions. This may lead to the natural tendency to study proof, often seen as the necessary leap from computational mathematics into the theoretical, independently by controlling for any confounding variables that may appear in the form of difficult non-proof concepts. Of course, a non-proof concept is typically the subject of proof construction, so researchers may seek students who "know" the non-proof concept that will be the subject of the proof task in their study (e.g., Healy & Hoyles, 2000; Weber, 2001). What is implicitly assumed in this type of study design is that any difficulties a student faces with proof construction will mostly be a result of student difficulty with proof only; investigating the nature of that difficulty can further illuminate the reasons why students may struggle with proof as a whole.

This method, though reasonable, may still pose an issue if the researchers involved in the study want to make the claim that a student "understands" a mathematical concept or "possesses" the ability to construct proofs. This type of claim is typically made by using the criteria that a student has arrived at a "correct" solution, can recall definitions, or has encountered content previously as a proxy for student "knowing", "learning", or "understanding", even in the absence of operationalized definitions of these terms. For whom a solution is deemed correct is seldom explicit, leading to ambiguities surrounding any attempts to assess and make conclusions about a student's *concept* of proof, their *concept* of non-proof mathematical content, and the bases on which conclusions have been made concerning the state of these *concept*s.

Such ambiguity is to be expected when researchers approach "knowledge," for example, from purely behaviorist, psychological, or sociocultural perspectives. As researchers, we are not able to see into the minds of students, though using objectified cognitivist terms obscures this fact, and it is necessary to rely on what is uttered by students to get a glimpse into how they interact with mathematics. In this research, I will incorporate analytical tools related to studying concept formation and structure that have been developed from the cognitivist perspective, but I will attempt to eliminate some of the ambiguity surrounding concepts in this perspective by strictly adhering to the discursive definition of *concept* as defined by Sfard (2008). I may also appeal to ideas from a behaviorist perspective to discuss changes in student activity, and ideas from the sociocultural perspective to capture the influence of societal interactions and external demands on students' goals for their mathematical behavior. However, I do not envision this research as a pragmatic mixture of these perspectives, but rather, I reposition student difficulty with proof and student difficulty with the elements of function as two instances of the same problem from a communicational participationist standpoint: potentially incommensurable discourses between discursants.

In this brief overview of student struggles with proof and function, I have introduced the epistemological and ontological issues that accompany previous attempts to study student thinking. In response to these issues, I have included terminology from Sfard's commognitive framework to prime the reader for the view of mathematics as discourse that appears necessary for the rigorous unobjectified study of students' mathematics (Sfard, 2008). While a possible outcome of the study that I undertake is insight into ways students may meet the expectations of introductory proof course instructors without necessarily acclimating to mathematical community norms about injective functions, surjective functions, nor proof, the primary focus of the research proposed herein is to gain insight into ways students organize elements of mathematical discourse to make sense of injective and surjective functions, interpret the demands of associated proof tasks, and ultimately engage in proof construction in response to these tasks.

**Research Questions**

1. In what ways do undergraduate mathematics students interpret and make sense of injection, surjection, and proof instruction in transition to proof courses?
    - How do students think about injections and surjections in the context of function before participating in transition to proof courses as evidenced by their initial discursive maps?
    - What changes do students make to their discursive maps pertaining to proof, injections, and surjections as a result of lecture, assignments, and reading of the course textbook as well as the use of other student-selected resources over time?

- How do students draw upon their injection, surjection, function, and proof discourses to interpret proof task situations involving injections and surjections?
- What routines do students draw upon to construct proofs in response to task situations?

2. In what ways do instructors interpret and make sense of proof, injective functions, and surjective functions?
    - How do instructors think about injections and surjections within the context of function as evidenced by their discursive maps?
    - How do instructors expect students' injection and surjection discourses to be organized within the context of their function discursive maps prior to the transition to proof course? Including proof discourse, how do instructors expect students' discourses to be organized by the end of the transition to proof course?
    - What discourses do instructors expect students to draw upon to interpret task situations involving injections and surjections in order to engage in proof construction as evidenced by their expected student trees?

3. How do the discourses of students differ from the discourses of their instructors?
    - How do instructor's expectations of students' injection, surjection, function, and proof discourses prior to and at the end of the course align with students' discourses on injections, surjections, function, and proof as evidenced by their initial and final discursive maps?
    - For each task situation, how do instructor's expectations of the discourses that students will draw upon to interpret proof task situations involving injective and surjective functions and engage in proof construction align with the discourses that students' draw upon as evidenced by their evoked trees and discursive maps?
    - How do instructor's expectations of students' function, injection, surjection, and proof discourses based on their reading of students' proof constructions align with students' discourses on function, injections, surjections, and proof as evidenced by their discursive maps?
    - How do instructor's narratives about students' proof and function discourses after reading students' proof constructions align with their narratives on students' proof and function discourses based on students' discursive maps and evoked trees?

**Chapter 2: Theoretical Framework & Literature Review**

**Sfard's Theory of Commognition as a Theoretical Framework**

Sfard (2008) describes commognition as a "term that encompasses **thinking** (individual cognition) and (interpersonal) communicating; as a combination of the words **communication** and cognition, it stresses the fact that these two processes are different (intrapersonal and interpersonal) manifestations of the same phenomenon" (p. 296, emphasis in original). A difficulty within mathematics education research is that much of the phenomena of interest (e.g. student thinking) are thought to be mental or abstract in nature. Furthermore, the mathematical objects of interest are themselves often abstract. Take, as a general example, what is implied by the phrase "has a learning disability," or in the subject of this research, what it means to be an injective function. These words do not describe physical and tangible phenomena, rather, they are the result of human *reification* and *objectification* of abstract notions over time subject to persistent ontological collapse. Sfard (2008) defines *reification* as "the replacement of talk about processes with talk about objects" (p. 301) and notes that this may involve the use of a new word, term, or symbolic artifact to discuss that process. Related to reification is the act of *saming,* defined as assigning a singular word, term, or symbolic artifact to multiple objects not previously considered the same, but that are interchangeable under specific circumstances (Sfard, 1998, p. 302). After reification the resulting term, word, or symbol may then result in the process' *objectification,* which Sfard describes as the use of a symbolic artifact or noun as though it "signified an extradiscursive, self-sustained entity (object), independent of human agency" (p. 300).

As a basic example, the *numeral 7* is often taken as a self-sustained entity by adults in the population of focus. It is no longer remembered that the notion of *7* needed to be built within the adult when they were a child, typically starting with the use of a number song to determine a quantity of physical or visually perceptible items. Ending the number song when there are no longer physical objects to assign a number word, a child in this population may announce that there are *seven* items. Through the repeated exposure to a series of objects that, when a child is instructed to *count*, do not require the use of number words after *seven* in the number song, a specific quantity of objects comes to be associated with the number word *seven* within the child. In an act of subitizing (Beckmann, 2013), the child may come to perceive a quantity of items as *seven* without appealing to their number song or other counting strategies. *7* is then introduced as a visual symbol to be associated and used interchangeably with the number word *seven*. Yet, it is taken for granted that a written phrase like *7 apples* is a simple description of a physical truth, rather than the result of several lengthy discursive maneuvers over time. What started as multiple physical tangible objects to the child becomes a permanent reference to be associated with the numeral 7. This occurs even though 7 is not a *primary object* (p-object), but instead a discursive object (d-object). *P-objects* are inclusive of what can be seen, touched, or only heard. *P-objects* must be tangible, perceptually accessible, and not brought into being by human communication; even without a signifier, they must exist independent of human discourse (Sfard, 2008, p. 169).

**REDACTED IMAGE**

Figure 1: Primary object (p-object)

Although I am including an image of the object in Figure 1, it is one with which many people come to be familiar and have observed physically in their lifetime. In the discursive context, the object in Figure 1 depicts a primary object (p-object). A *d-object* can be a primary object or a compound discursive object that "arises by according a noun or pronoun to extant objects in the process of **saming, encapsulation or reification**" (p. 298, emphasis in original). In this sense, a d-object can arise by assigning a single p-object to a set of p-objects; alternatively, a d-object can also arise by assigning a single p-object to a set of previously constructed d-objects. While it is still unnamed, the object in Figure 1 can only be referred to as a p-object, not a d-object. I note that the noun typically associated with the object in Figure 1 is *Gala apple* in English. As a result of the object in Figure 1 now being signified by *Gala apple* and thus being an object of communication, the object in Figure 1 can be referred to as a d-object. *Gala apple* does not exist independently, but rather, the noun *Gala apple* is the result of our culture's fluid and constantly revised communicational system. Without our communicational system, *Gala apple* does not exist, while the object in Figure 1 remains tangible.

The term *Gala apple* is not naturally endowed with the series of characteristics that we require before we call an object *Gala apple*, rather, the collection of specific characteristics "possessed" by the observed primary object was assigned the term *Gala apple* by a prior culture. There is only one p-object that is called *Gala apple*, and therefore, *Gala apple* can be called an *atomic* or *simple d-object*. A *simple d-object* is defined as any d-object that results from proper naming of a p-object (e.g., *Gala apple*) where a term is exclusively associated with a particular p-object (Sfard, 2008). The simple d-object in this pair now serves as a *signifier* that can be used to effectively communicate about one specific p-object (e.g. the p-object in Figure 1). The p-object in this example serves as the only realization for the signifier *Gala apple*. Note, the signifier *apple* would not be referred to as a simple d-object due to numerous realizations in the form of different pronouns for, or "types" of, *apple* that are typically brought to mind when the term *apple* is used. Instead, *apple* can be called a *compound d-object*.

A *compound* d-object results from again assigning "a noun or pronoun to extant objects (either discursive or primary) by **saming, encapsulating, or reifying**" (p. 296, emphasis in original). Any p-object that is labeled *Gala apple* can also be called *apple*, and all characteristics that are used to identify a p-object as *apple* are contained in the characteristics we use to label an object *Gala apple*, but there exist other p-objects that we label with d-object *apple*, but that we do not label *Gala apple*. In commognition, *apple* arises from *saming*. This saming results from the identification of a specific subset of characteristics that are shared by a set of p-objects, where each p-object may have been assigned a simple d-object such as *Gala apple*, *Macintosh apple, Green apple*, etc. and assigns any p-object with those shared characteristics the noun *apple*.

*Encapsulation*, defined as the "act of assigning a noun or pronoun (**signifier**) to a specific set of extant **primary** or **discursive objects**, so that some of the stories of this set that have, so far, been told in plural may now be told in singular" (p. 298, emphasis in original) can be exemplified by the d-object *fruit cup.* If an object is called *fruit cup,* stories about the *fruit cup* are specific to the elements of the *fruit cup* when taken together. If *fruit cup* is referred to as *expensive,* the intent of this communication is about *fruit cup* as a whole. It does not apply to the individual items that contribute to *fruit cup*. I emphasize that objects to be encapsulated or samed do not need to be primary objects nor simple d-objects. Several compound d-objects can again be encapsulated or samed into a new compound d-object.

The term *concrete object* is reserved for the specific cases where the objects to be encapsulated or samed are either p-objects or the d-objects which arise through the saming and encapsulation of those p-objects. For example, consider the signifier *fruit*. *Fruit* is not only a compound d-object, but it can further be described as a *concrete object* that arises from the encapsulation and saming of d-objects such as *kiwi*, *pineapple*, *grape*, etc. and the p-objects associated with each of these nouns. *Kiwi*, *pineapple*, and *grape* are individually concrete objects resulting from the saming and encapsulation of p-objects (Sfard, 2008). Concrete objects can be repeatedly broken down into other concrete objects and this process terminates upon reaching a simple d-object, where the sole realization of this simple d-object exists in the form of a unique p-object. In short, concrete objects must be perceptually accessible.

The complexity of human language is striking due to its ability to communicate ideas about objects that have no pure basis in the physical world. Any object that is not a concrete object, that is, an object that arises in whole or in part from the encapsulation and saming of purely discursive objects, is referred to as an *abstract object*. For example, *7* in the phrase *7 apples* is an *abstract object*. Recall the previous lengthy discussion on how an adult comes to a notion of *7* over a period of time during their childhood; I sought to demonstrate that *7* does not exist independent of discourse. This is the opposite of the concrete existence of the p-object in Figure 1. The language of choice for which the signifiers for this p-object are selected is irrelevant. The p-object in Figure 1 does not cease to exist if it does not have a proper noun associated with it. While *Gala apple* is a d-object resulting from reification, the p-object that is realized upon reading or hearing *Gala apple* is typically one similar to what is depicted in Figure 1. This is true for all p-objects that are referred to as *apple,* and the term *apple* and all associated p-objects are therefore concrete objects. An attempt at writing out an abstract object's *realization tree,* defined as a multi-layered, hierarchical, written representation of how a term or symbolic artifact can be deconstructed into its successive object components, cannot terminate by arriving at solely atomic d-objects (Sfard, 2008). I invite the reader to attempt this.

While concrete objects like *apple* and *Gala apple* are commonly taken for granted and regarded as indisputable preexisting "things in the world" (Sfard, 2008), difficulties arise when attempting to ask for an *apple* from a non-English speaker or, when speaking to someone for which *apple* describes a completely different p-object than one would associate with it in the United States. The notion of *apple* and *Gala apple* are only relevant in a specific c*ommunity of discourse*, or the collection of individuals who can participate in a discourse. The community of discourse in this example may be English speakers, or if I am being cautious, English speakers who grew up speaking English in the United States. The issue I describe regarding how an English speaker and non-English speaker may hold different realizations for the signifier *apple* is described as having *incommensurable discourses.* Sfard describes discourses as incommensurable if the discourses "differ in their use of words and **mediators** or in their **routines**; incommensurable discourses may allow for the endorsement of seemingly contradictory **narratives**" (p. 299, emphasis in original).

Extended participation in a common discourse community often leads to communicationally advantageous *ontological collapse*. Ontological collapse is defined by Sfard (2008) as the "phenomenon of taking all the objects that are being talked about – discursive (d-objects) and extradiscursive (p-objects) – as belonging to the same ontological category of 'things in the world' that preexist discourse, with their mutual relations similarly "objective" and mind-independent" (Sfard, 2008, p. 300). This lack of awareness of language as constructed can be seen in, for example, cultural replacement of the notion of a p-object that has been termed "Gala apple" with a statement that the object in Figure 1 *is a* Gala apple; a language independent and pre-existing thing-in-the-world. Attempting to discuss the creation of all d-objects from only their primary objects is, as I believe I have demonstrated, a lengthy and often unnecessary task when in the company of fellow members of a discourse community. Ontological collapse enables efficient communication about increasingly complex ideas within that community. However, resisting ontological collapse and distinguishing between objects in the manner that I have outlined is helpful when it is necessary to delineate between tangible physical objects, observable actions, and purely mental objects. These distinctions are necessary for my study.

For the sake of specificity in discussing the literature related to this dissertation, objects need to be deconstructed even further. I remind the reader that Sfard defines *realization* as a procedure that pairs a *signifier* with another object (Sfard, 2008) and defines *signifier*s as "words or symbols that function as nouns in utterances of discourse participants" (Sfard, 2008). As a note, the distinction between these two terms is often relative, and realizations can serve as signifiers for another set of realizations. Furthermore, realizations are limited to what is perceptually accessible. This is what thwarts the task of completing a realization tree for the abstract d-object *7*. The leaves of such a tree can only be a never-ending list of purely discursive objects.

**REDACTED IMAGE**

Figure 2: Realization modalities as listed by Sfard, represented in the format of a partial realization tree. As a note, these modalities do not represent the full spectrum of communicational types. Reprinted from Sfard, A. (2008). *Thinking as communicating: human development, the growth of discourses, and mathematizing*. Cambridge University Press.

Now that signifiers, realizations, p-objects, simple (atomic) d-objects, compound d-objects, and the abstract or concrete nature of d-objects have been carefully defined, I note that I have only discussed how d-objects arise. I have not discussed how d-object is operationalized. Due to the prevalence of the term *object* within mathematics education research, it is important that I distinguish how I will use this term within my study. According to Sfard, "The (discursive) object signified by S (or simply object S) in a given discourse on S is the realization tree of S within this discourse" (Sfard, 2008, p. 166). This definition provides specificity, the possibility of symbolic representation, and an element of recursion. It indicates that the domain of discourse is relevant to the realization tree. Furthermore, realizations, and in turn realization trees, are unique to the person making them, so "uses" are also person-specific. This is especially helpful when defining mathematical objects.

*Mathematical objects* are abstract objects that are "complex hierarchical systems of partially exchangeable symbolic artifacts" (Sfard, 2008, p. 172). In turn, mathematical discourse can be defined as discourse about mathematical objects. *7* is an example of a mathematical object, and similarly, so are *3, 4,* and *+*. When operating within *natural numbers modulo four*, the symbolic artifacts *7* and *3* are exchangeable (i.e. interchangeable). However, this would not be true when operating within the *field of real numbers*. This is what necessitates the inclusion of the term *partially* in Sfard's description of the interchangeability of mathematical objects. Furthermore, grappling with the real numbers subsumes the discourse of integers, which itself subsumes the natural numbers, which in turn can be realized as specific numerals, where these numerals are symbolic artifacts borne of lengthy experiences with words in a number song. A hierarchy is established each time mathematical objects are samed, encapsulated, or reified to produce new mathematical objects. At each level of encapsulation, there is an increase in the complexity of the new mathematical object that is created, and there does not appear to be an end to the to the constant emergence of new mathematical objects as a result of this process. For this reason, Sfard describes mathematics as *autopoietic*, which is defined as "a system that contains the objects of talk along with the talk itself and that grows incessantly 'from inside' when new objects are added one after another" (Sfard, 2008, p. 129). In short, mathematical discourse is self-generating. Mathematical objects *are* realization trees.

If mathematics is built upon a foundation of abstract objects, how can it be consistently ensured that members of the mathematics discourse community have the same or compatible concept (i.e. a signifier together with its uses) for each mathematical object? For more foundational concepts such as *7*, *subtract,* or *add*, creation of tasks to examine a person's usage when confronted with these terms or situations that call for them may seem easy. But how can we accurately determine a person's concept of *isomorphism, ring,* or *integral*? This is decidedly more difficult than determining a person's uses of and for *snake, tree,* or *vehicle*. *Isomorphism, ring,* and *integral* cannot be pointed out while a child walks to school. We would rely entirely on what someone communicates about these abstract objects if we were to try to gain insight into their concept of these objects.

Much of mathematics education research revolves around the thinking of students. However, thinking within an individual is hidden to anyone outside of the individual. When a researcher describes the thinking of a participant within the population of focus, this is at worst mere interpretation of observable activity and at best a word-for-word replication of what a participant communicates that they are thinking. It is therefore sensible to examine thinking and communication in tandem. Sfard's (2008) theory of commognition is built upon this idea. Under the participationist view that "human development is the process of individualization of historically established, collectively implemented forms of life," the central tenet of commognition is that *thinking* is a form of communication, and this communication can be both interpersonal and individual (p. 262). While this communication need not be verbal, the recursive self-referential nature of human language allows for further complexity of our activities and our communication about those activities over time. This communication is often aided by ontological collapse.

Under this theory, Sfard positions mathematics as *discourse*, or a "special type of communication made distinct by its repertoire of admissible actions and the way these actions are paired with re-actions" (p. 297). Every distinct discourse constitutes a community (e.g. academic discourse; social justice discourse) and discourses in language are "defined by vocabulary, endorsed narratives, visual mediators, and routines" (Sfard, 2008). Mathematical discourse can be distinguished from other discourses based on any discursant's ability, in theory, to endorse or reject narratives; and on the use of vocabulary that is distinctly mathematical in nature (i.e. mathematical objects.) The endorsing authority within mathematics discourse does not lie within specific individuals, but rather, any member of the discourse can substantiate a narrative by appealing to some combination of endorsed narratives. In order to maintain the integrity of Sfard's framing of mathematics as discourse, I have included Sfard's brief descriptions of the defining features of mathematical discourse in whole below with the caveat that some features are exclusive to specific communication communities:

1. *Word use.* One of the distinctive characteristics of discourses is the keywords they use. In mathematics, these are mainly, although not exclusively, the words that signify quantities and shapes. Whereas many number-related words may appear in non-specialized, colloquial discourses, mathematical discourses as practiced in schools or in academia dictate their own, more disciplined uses of words. Word use is an all-important matter because, being tantamount to what others call "word meaning," it is responsible for what the user is able to say about (and thus to see in) the world.
2. *Visual mediators* are visible objects that are operated upon as a part of the process of communication. While colloquial discourses are usually mediated by images of material things existing independent of the discourse, scientific and mathematical discourses often involve symbolic artifacts, created specifically for the sake of this particular form of communication; think, for example, about scientific inscriptions or mathematical algebraic notation. Communication-related operation on visual mediators would often become automated and embodied. Think, for example, about the procedures of scanning the mediator with one's eyes in a well-defined way. With some experience, this procedure would be remembered, activated, and implemented in direct response to certain discursive prompts, as opposed to implementation that requires deliberate decisions and the explicit recall of a verbal prescription for these operations.
3. *Narrative* is any sequence of utterances framed as a description of objects, of relations between objects, or of processes with or by objects, that is subject to *endorsement* or rejection with the help of discourse-specific substantiation procedures. Endorsed narratives are often labeled as true. Terms and criteria of endorsement may vary considerably from discourse to discourse, and more often than not, the issues of power relations between interlocutors may in fact play a considerable role. This is certainly true about social sciences and humanistic narratives such as history or sociological theories. Mathematical discourse is conceived as one that should be impervious to any considerations other than purely deductive relations between narratives. In the case of scholarly mathematical discourse, the consensually endorsed narratives are known as mathematical theories, and this includes such discursive constructs as definitions, proofs, and theorems.
4. Routines are repetitive patterns characteristic of the given discourse. Specifically, mathematical regularities can be noticed whether one is watching the use of mathematical words and mediators or following the process of creating and substantiating narratives about numbers or geometrical shapes. In fact, such repetitive patterns can be seen in almost any aspect of mathematical discourses: in mathematical forms of categorizing, in mathematical modes of attending to the environment, in ways of viewing situations as "the same" or different, which is crucial for the interlocutors' ability to apply mathematical discourse whenever appropriate – and the list is still long. (Sfard, 2008, pp. 134-135)

Sfard's description of narrative is important as endorsement or rejection is specific to the discourse. Mathematical discourse is composed of numerous discourses; what can be endorsed in one discourse may be rejected in another. For this reason, Sfard's differentiation between mathematical discourse and scholarly mathematical discourse is appropriate. Though narratives are of importance, what I have found most interesting in Sfard's description of mathematics as discourse is the discussion of routines. If I define concept as a signifier together with its uses, then I see routines as the building blocks of concept, and a collection of all routines paired with a signifier and its realizations as forming a person's concept of a mathematical object. Word use, visual mediators, and narrative are themselves elements of routine. Routines can also be present within larger routines.

### Routines.

Of particular importance to my research and data analysis is the notion of routines when viewing mathematics as discourse. *Routines* are defined as a "set of metarules that describe a repetitive discursive action" (Sfard, 2008, p. 208) and can be broken down into three features: applicability conditions, course of action, and closing conditions. *Applicability conditions* dictate "when" a routine should be used and are the circumstances a discursant seeks to be met before enacting a particular routine. A *course of action,* or *procedure*, is the "how" of the routine that is enacted under the specific set of applicability conditions. *Closing conditions* are the circumstances a discursant seeks before the routine can be terminated. This *closure* defines another "when," specifically in reference to when a routine can be considered completed to the satisfaction of the discursant.

### *Types of routines.*

According to Sfard, the term *endorsable* "signals that the narrative can be endorsed or rejected according to well defined rules of the given mathematical discourse" (p. 224). Examples of the most basic endorsed narratives in mathematics are axioms and definitions. When these axioms and definitions are used according to well-defined rules within the discourse to derive yet another narrative, the derived endorsed narrative is typically referred to as a theorem (p. 224). The collection of several such endorsed narratives forms a system of mathematical theories (p. 224). During the course of a routine, discursive and practical actions may be present in addition to the manipulation of narratives. While discursive actions are those that are communicational in nature, *practical actions* result in a physical change in objects or environment (p. 236). Now that definitions have been provided for these terms, I can expand on routine's different types: explorations, deeds, and rituals.

*Explorations* are routines that terminate when "an endorsable narrative is produced or substantiated" thus contributing to mathematical theory. There are three subtypes of exploration routines (p. 224): (1) *Construction* explorations are discursive processes that result in new endorsable narratives; (2) *substantiation* explorations are the actions that aid in the decision to endorse a previously constructed narrative; and (3) *recall* explorations are the "process one performs to be able to summon a narrative that was endorsed in the past" (p. 225). These three types of exploration are commonly interwoven in the process of endorsing a narrative. Of note, substantiations and recall rely on previously endorsed narratives within the discourse. The opportunity for creativity in mathematics lies within construction explorations.

*Deeds* are defined by Sfard as "a set of rules for a patterned sequence of actions that, unlike explorations, produce or change objects, not just narratives" (p. 236). This type of routine terminates when the participant in the discourse of mathematics (i.e. the discursant), hereafter referred to as a *mathematist* (Sfard, 2008), believes they have transformed the environment or object in fulfilment of a practical task. There is no goal of producing a narrative about the object, but rather, the transformation of the object is the goal. Note, the same set of instructions provided to two mathematists may result in the enaction of an exploration by one mathematist, but a deed by the other mathematist.

*Rituals* are defined by Sfard as, "sequences of discursive actions whose primary goal (closing conditions) is neither the production of an endorsed narrative nor a change in objects, but creating and sustaining a bond with other people" (p. 241). Such routines place emphasis on social approval and can be seen in students who are preoccupied with "getting to the answer the teacher will mark correct." These routines can entail performance for the perceived authority figure and can be highly situational. However, ritual can also be seen in an expert's normed use of specialized words to refer to mathematical objects. By adherence to ritual, the expert can assure that they are communicating what they intend to another expert, thus sustaining a bond to the mathematicians that came before them.

Routines can be mathematist-specific, although multiple mathematists may apply the same applicability conditions and closing conditions to a similar procedure. Within a *task situation*, defined as any situation in which in a person believes they are required to do something (Sfard, 2008; Lavie, Steiner, & Sfard, 2018), applicability conditions may determine the *prompts* that incline a mathematist to enact a particular routine. According to Sfard, prompts are

> the elements of situations whose presence increases the likelihood of a routine's performance. The prompts, which may be verbal or environmental and can be provided by others or self-presented, occur in clusters. Only some combinations of prompts suffice to spur routine performance. When a person seems to be following a given course of action restrictively, showing a greater dependence on situational clues than would be proper in the eyes of an experienced mathematist, we say that her discourse (or just the given course of action) is *situated*…Some of the prompts associated with a given course of action function as *discourse framers*, that is, as factors that activate certain discourses while inhibiting some others. (Sfard, 2008, p. 209, emphasis in original)

Note, similar procedures may also be enacted by two mathematists despite applying different applicability conditions and having distinct conditions of closure. For this reason, a student's execution of the procedure expected by the instructor should not be taken as evidence that the student would enact the same procedure in circumstances considered appropriate for a more experienced expert mathematist. However, this is precisely what is occasionally assumed in mathematics classrooms and mathematics education research.

**Literature Review**

      Upon analysis of common methods used in mathematics education research, learning and understanding are routinely assessed using the implicit criteria of adherence to *norms* (e.g. Breidenbach et. al.,1992; Dubinsky & Wilson, 2013; Wilson, 2007); this is evidenced by use of student work alignment with normed constructions, "correct" answers, and student recall of definitions in the methods and subsequent analysis as indicators, or lack thereof, of student progress. Norms can be defined as *metadiscursive* rules, hereafter referred to as *metarules*, that are widely endorsed and enacted with the discourse community (Sfard, 2008, p. 204). For the purposes of this research, the discourse community is the mathematical community. Furthermore, this community includes perceived experts, mathematicians, mathematics professors, mathematics educators, as well as student participants in the mathematics discourse. When specificity is necessary, student participants will be referred to as student mathematists, while mathematicians and mathematics professors will be referred to as expert mathematists. *Metarules* "define patterns in the activity of the discursants trying to produce and substantiate object-level narratives" (p. 201) and *object-level rules* are "narratives about regularities in the behavior of objects of the discourse" (p. 201). Therefore, when I use the term norms, I am referring to the patterns of activity widely accepted by expert mathematists when trying to produce and substantiate narratives about mathematical objects. The term *classroom norms* will be used to refer to the patterns of activity endorsed by the professor of the mathematics course in agreement with or to be expected of the student mathematists within the classroom.

      When a researcher is attempting to examine student mathematists' adherence to norms, there is a tendency to rely on the written or spoken utterances of the student mathematist for evidence of adherence to those norms. While I have stated that I agree that written and spoken utterances can provide helpful information about student thinking and concept, any attempt to study student adherence to norms further requires an explication of the conditions that each student mathematist uses to govern their activity. In response to the lack of this sort of explication in the literature, I note that the metarules of the expert mathematical community (norms) often differ from that of the student mathematist (e.g., Erlwanger, 1973). If this difference is present, it would be reasonable to expect that misalignment in metarules on the part of the student when compared to experts would become immediately obvious to the researcher, but this expectation relies on the student enacting a written or spoken *routine* that lies in conflict with the routines that would typically be employed. If the student mathematist enacts a routine that is similar to or exactly the routine that is expected, the researcher and the course instructor may accept this routine as evidence of adherence to the norms of the mathematical community. In turn, it may then be assumed that the student has "learned" or "understands."

      As mentioned, this latter point, while interesting, is not the focus of this dissertation, although the commognitive view does operationalize the term *learn* by defining learning as a change in discourse. Meditations on what it means to learn or understand have been covered in depth elsewhere (e.g. Sfard, 2001, 2008; Barwell, 2007), and I find myself in alignment with Sfard's argument that the common tools utilized to assess learning and understanding are more so assessing discursive competence (Sfard, 2008). The former point, that the exact or similar enactment of routines by students may occur outside of normed mathematical discourse, is an interesting possibility that may be encountered when examining how students organize and draw on elements of their mathematical discourse. For this reason, I will use the term "sense-making" in place of learning and understanding to account for all student mathematical discourses, and

not only those discourses that would be deemed satisfactory to a large portion of typical expert mathematists.

**Research on sense-making.**

Sense-making, as described by Schoenfeld (1991) in his seminal work, parses out the multiplicity of ways in which learning and understanding are used to describe how students process and respond to mathematical information while additionally allowing for doubt to be placed on both performance and behavior as reliable indicators of adherence to the accepted norms of the mathematical community. The types of sense-making as alluded to by Schoenfeld are described as follows, and these categories are not necessarily disjoint.

*Imparting meaning to symbols.*

In this view of sense-making, Schoenfeld refers to student ability to translate between contexts into mathematical language. He states, "…the act of transforming the story into the arithmetic operation...is an act of sense-making; an act of imparting meaning to mathematical symbols." (Schoenfeld, 1991, p. 322) Once a student arrives at symbols, they can then perform symbolic manipulations based on the formal system in which they are expected to operate. Of such systems, Schoenfeld states, "formal systems are not about anything. Formal systems consist of sets of symbols and rules for manipulating them. As long as you play by the rules, the results are valid within the system" (Schoenfeld, 1991, p. 311). In this view of sense-making, students can possibly "game" the process through the use of keywords. For example, a student in an introductory proof class may associate being asked to *find the limit* of a function using the delta-epsilon definition with setting up an inequality due to watching a series of examples where this exact step is demonstrated by their teacher. This does not, however, indicate that a student connects this process to the definition, nor that they have assigned a meaning to what this process produces that aligns with mathematics community norms. In an elementary case, this can be exemplified by a student who sees the word *more* in a word problem and immediately places a plus symbol in between the number symbol for number words they can identify in the problem, without regard for additional information in the problem that may suggest that solely performing addition is inappropriate. Weber (2001) referred to this type of sense-making as syntactic knowledge or syntactic understanding.

*Behavior.*

On behavior as a type of sense making, Schoenfeld states:

> behavior is sense making of the deepest kind. In the context of schooling, such behavior represents the construction of a set of behaviors that results in praise for good performance, minimal conflict, fitting in socially, and so forth. What could be more sensible than that? The problem then, is that the same behavior that is sensible in one context (schooling as an institution) may violate the protocols of sense making in another (the culture of mathematics and mathematicians). (Schoenfeld, 1991, p. 340)

In this view of sense-making, Schoenfeld highlights the impact that an external response, whether positive or negative, can have on students' mathematics. Cobb (2000a) further illustrates this impact in his commentary on classroom-based research that focuses on students' mathematical learning and development:

> The classroom microculture in which students participate influences profoundly the goals they attempt to achieve, their understanding of what counts as an acceptable mathematical explanation, and indeed, their general beliefs about what it means to know and do mathematics in school. Consequently, it is essential to document the microculture established by the classroom community even if the primary concern is to formulate psychological models of the processes by which students transform their mathematical ability. (Cobb, 2000a, p. 311)

Students may change solution approaches based on what they believe their teacher wants to see from them, and will discard approaches, whether they are deemed correct or not, if a positive response (e.g. earning a "B" grade) is not achieved, or if a perceived negative response in the form of a correction by the instructor is received. For a student to ensure that a produced solution is in line with teacher expectations, it is possible, and even common in proof, to see students resort to the memorization and mimicking of processes they see their teachers engage in, without regard for the meaning behind what they have produced (Moore, 1994). For example, a student may leave an introductory Calculus course with the ability to compute limits using limit rules that are standard for high school Calculus and that have been labeled correct by their Calculus teacher, but upon entering a transition to proof course, they may seemingly begin to appeal to the delta-epsilon definition only after several demonstrations by their teacher where the delta-epsilon definition is invoked. This student, though able to compute a limit using Calculus rules and able to prove that a limit exists using a delta-epsilon argument, may nevertheless be unable to discuss the relationship between these rules and the delta-epsilon definition of a limit.

### *Pattern identification and rule creation.*

In this view of sense-making, a student may identify common features in similar sets of problems and attempt to determine rules for solving problems based on these features. If a teacher is prone to assigning problem sets that give students the opportunity to practice newly taught material immediately after the material is introduced, a student may determine how to arrive at a "correct" solution on a problem set and then consistently apply this technique to all problems in that set. For example, Schoenfeld refers to a group of elementary school students who were provided with several change unknown subtraction problems in the form a - _ = b after first learning about subtraction. The problems did not deviate from this format. A second set of problems included a variety of start unknown, change unknown, and result unknown subtraction problems. Of the first problem set, Schoenfeld states, "Such problem sets have an interesting and unintended property. Once you figure out how to solve the first problem (in this case, subtract the smaller number from the larger), that method works on all the problems" (Schoenfeld, 1991, p. 321). However, on the second problem set, students still appealed to the rule set created for the first problem set, thus leading to several "incorrect" and nonsensical answers. Again, I note that a student being able to determine how to arrive at a "correct" answer does not imply that their process for arriving at a correct answer was itself logical in normed mathematical discourse. This view of sense-making can be seen as connected to behavior as sense-making, but in this view, a student may create and adhere to rules without necessarily receiving confirming input from their teacher. As can be seen in Erlwanger (1973), rule sets created by students that result in "incorrect" solutions only some of the time will not necessarily lead a student to discard their rule set. These idiosyncratic rule sets can be referred to as *metarules* and in this instance, they refer to the student's actions and patterns in their activity when attempting to substantiate *object-level* narratives (Sfard, 2008, p. 201).

### *Adherence to established mathematical culture.*

In this view of sense making, students are generally regarded as having adopted interpretations of content that are consistent with the established mathematical culture (norms). These students are not viewed as merely guessing or mimicking what is expected them, but as having made a connection to their existing and (presumably "correct") mathematical knowledge (discourse.) Schoenfeld alludes to the absence of this type of sense-making in some students by describing a situation with high school students that is at the crux of why it is important to study sense-making in proof

> To generate the context for a discussion of high school students' understanding of geometry, I gave the students a proof problem to solve. The problem...was written on the chalkboard. Without hesitation, the students produced a complete and correct solution to the problem. They dictated the answer, and I wrote the proof on the board. The whole process took less than three minutes from start to finish. (Schoenfeld, 1991, p. 317)

Schoenfeld then goes on to ask his students to draw the accompanying construction for the proof they have just written. On this second task, four solutions were proposed, but only one student was able to determine which construction was correct. Schoenfeld suggests that this failure to connect proofs to their actual "meaning" was common, "I was surprised to see it emerge with such talented and advanced students, but there is ample documentation...that high school and college students see very little connection between the world of deductive Euclidean geometry and closely related construction problems that they are asked to solve" (Schoenfeld, 1991, p. 319). In this case, perhaps as a consequence of behavioral sense-making, Schoenfeld surmises that students memorized several proofs, any of which could appear on an end of year exam, and the reproduction of such a proof would earn the students significant points on their exam. While all of these students could produce the proof, students were unable to produce the construction (referring to the object described within their proof) by exploiting the necessary details in the proofs they had just written. Although students were "correct" in their proof, their *concept* of what they had just described did not allow them to make any further inferences in response to Schoenfeld's subsequent task.

Of Schoenfeld's identified types of sense-making, adherence to established mathematical culture (i.e. alignment with the concept image common to mathematicians in cognitivist terms) appears to be what is desired by most instructors for their students' mathematical development, with the added necessity of students appropriately imparting meaning to symbols. However, behavior, pattern identification, and rule creation by students may allow students to appear to be adhering to established mathematical culture when this may not accurate. As a student advances in the mathematics major, it would likely be increasingly difficult for that student's idiosyncratic rule sets to suffice. However, if a student possesses a highly sophisticated set of deeds and ritualized routines, this student's lack of normed mathematical discourse will likely go unnoticed.

Some examples of the various aims of studies on student sense-making are to examine the reasoning and choices students make in problem solving, the creation of an image of student's assumed "knowledge" structure and when it is accessed, as well as the changes that occur to "knowledge" structure over time (e.g., Adiredja, 2021; Plaxco & Wawro, 2015). Used more flexibly than terms such as learning, understanding, and knowledge in mathematics education research and unlikely to ignite ideological arguments over its colloquial use, the study of sense-making allows the researcher to take all relevant features of the environment into the analysis of students' mathematics.

Focusing on sense-making can provide insight into misalignment between written, enacted, and attained curriculum and offers teachers the ability to make changes to instruction that target various modes of sense-making. By examining the thinking of students and what influences that thinking, suggestions can be made concerning how the teaching of proof can be improved. Previous findings in the realms of beliefs, strategies of success/failure, and conceptions can provide the context for the "what and why" of student proof construction attempts, while a focus on sense-making can illuminate the "who, what, when, where, why, and how" of students' mathematics. Schoenfeld's analysis of sense-making raises important questions: (1) Do students know when they are not adhering to established mathematical culture? (2) If students are aware of this, is it intentional?

I note that many students literally do not see sense in mathematics. Schoenfeld (1991) states, "They trust their teachers; the scholastic setting of the problems coerces them into school-like behavior; they expect problems to be reasonable, to have answers derivable from the data given in them, and so on" (Schoenfeld, 1991, p. 320) This expectation leads students to abandon critical thought and the natural tendency to question what is presented to them. Students do not see themselves as having ownership in mathematics. If the teacher, the authority figure in the classroom, has presented the problem, then it is assumed that there must be a solution, which, perhaps, only the teacher will fully understand. Schoenfeld refers to this as being in a state of "non-reason" or the suspension of sense-making and defines suspension of sense-making as, "suspending the requirement that the problem statements make sense" (Schoenfeld, 1991, p. 316). He goes on to suggest that this "nonreasoning" is a result of the nature of schooling itself. He states

> The basic thrust of the argument is that classrooms (and more broadly, schools) are cultural milieux in which everyday activities and practices define and give meaning to the subject matter within them; culturally transmitted meaning-what students come to understand about mathematics as a result of their experiences within it-may or may not correspond to the intended meaning. (Schoenfeld, 1991, p. 320)

Student disposition toward a state of non-reasoning can have far reaching consequences for curriculum.

This difference in intended meaning and the meaning arrived at by students is well described in the context of curriculum by Kilpatrick (2009); curriculum can be seen as divided into three occasionally overlapping phases: the intended (elsewhere referred to as written), the implemented (elsewhere referred to as the enacted), and the attained or realized phase. The intended curriculum reflects the administrator's or curriculum writer's point of view, the implemented curriculum reflects the teacher's intentions and point of view, and the attained curriculum reflects the student's point of view (Kilpatrick, 2009). When there is harmony among the three phases, the expectation is that students attain and conceive of content precisely as intended by those responsible for the curriculum. The teacher's intentions for the content taught can be placed within the implemented or enacted phase, while the knowledge structures students possess at the conclusion of receiving instruction can be housed in the attained phase. The issue, therefore, is that students often attain "meanings" within curricula that are entirely different from what the instructor and curriculum writers intend. But what leads to this disharmony between the phases? If it is assumed that the teacher enacts the curriculum as it was written (intended), and yet the attained curriculum varies between students, then the disconnect between phases can possibly be attributed to differences at the level of individual student cognition and environmental influences on that cognition[3].

Schoenfeld poses a possible solution to the epistemological and ontological issues raised by the non-operationalized use of the terms learning and understanding within mathematics education research by uniting these uses under the name sense-making. However, Schoenfeld still makes extensive reference to "knowledge" despite using knowledge to describe states of the minds of students that cannot be seen by researchers. I re-emphasize that, at best, researchers can make inferences about student utterances (spoken or written) and draw conclusions based on those inferences. However, the pattern identification, rule creation, symbolic manipulation, behavior, and adherence to norms considered by Schoenfeld as the types of sense-making those students engage in can all be accounted for with the precise wording of Sfard's theory of commognition.

I connect Schoenfeld and Sfard's arguments in the following ways: (1) symbolic manipulation is accounted for in what Sfard calls deeds; (2) behavior and pattern identification are addressed within Sfard's routine which is inclusive of rituals; (3) pattern identification and rule creation can also be seen as student created narratives, object-level rules, and metarules; and (4) if object-level rules and metarules used by the student align with normed mathematical discourse, these would be considered endorsed narratives and norms under Sfard's theory; use of normed mathematical discourses is adherence to established mathematical culture. Schoenfeld's description of sense-making and Sfard's theory of commognition are therefore compatible, so I will operationalize sense-making as the collection of discursant-specific discourses around objects. In this study, these objects are mathematical in nature.

Sfard's (2008) theory of commognition is unique in that it allows the researcher to discuss student mathematics without the objectification of the students' mental state while also encompassing the ideas of the prominent theories and analytical tools used to study student mathematics. Commognition has been used to study, among many other topics, how 7th grade students develop their function discourse despite lack of a prior exposure to function (Nachlieli & Tabach, 2012), the use of dynamic geometry environments in increasing children's ability to identify shapes (Sinclair & Moss, 2012), changes in students' informal algebraic talk to formal algebraic talk (Caspi & Sfard, 2012), the equivalence of written and enacted curricula (Newton, 2012), and the impact of language on the learning of mathematics (Kim, Ferrini-Mundy, & Sfard, 2012). The unifying idea within commognition is that of the routine. To discuss the prior work that has been done in the areas of student proof, injections and surjections, and analyzing students' mathematics, I will occasionally make use of the terminology used by these researchers. When appropriate, I will operationalize objectified vocabulary using the terminology of commognition or clearly distinguish between commognitive usage and alternate usage of terms.

---

[3] I include embodied cognition (Nunez, Edwards, & Matos, 1999) within cognition. A discussion of my perspective in contrast with Sfard's view of embodiment as automation can be found in my methods.

**Research on proof and proving.**

Attempts to identify potential sources of student struggles within proof (e.g., Moore, 1994; Dreyfus, 1999) as well as the proposal of solutions to these struggles (e.g., Carrascal, 2015; Edwards & Ward, 2004) have been numerous over the past several decades. Selden & Selden (2008) theorized that students may reduce or eliminate difficulty with proof construction by entering post-high school studies after having some experience formulating simple proofs. However, upon examining beliefs of college students and mathematicians in the United States, Weber & Mejia-Ramos (2014) argue that students' unproductive beliefs about proof may have been acquired during such pre-collegiate proof experiences. Similarly, Coe and Ruthven (1994) found that college students in the United Kingdom who went through a proof-focused secondary curriculum relied on proof strategies that were largely routinized (used in a manner similar to the commognitive *ritual* as opposed to commognitive *routine*) and dependent on exposure. Furthermore, even when these students did construct proofs, these proofs were found to be restrictive (p. 52). Similar findings in the past 25 years would suggest that early introduction to proof is not currently solving the issue of student proof difficulty at the postsecondary level and may potentially exacerbate it.

When focusing on proof at the postsecondary level, the lines of inquiry in proof research can be grouped into four main categories: studies of student attempts and perceived success or failure, student and instructor beliefs, theories of proof, and the study of cognitive changes and sense-making. These categories are occasionally intertwined.

*Studies on success, failure, and attempts.*

Studies that seek to examine success in proof construction often involve participants who are perceived by the researcher as having the ability to engage in proof and attempt to illuminate the methods these participants use to construct a proof (e.g., Savic, 2015). The participants in these studies may be students, pre-service teachers, graduate students in mathematics, or possibly mathematicians. The methods of these participants may be self-reported (e.g., Duffin & Simpson, 2002). The designers of the study may then propose solutions to the problem of student proof construction difficulty by creating a framework or map that encapsulates the strategies of these participants (e.g., Zazkis, Weber, Mejia-Ramos, 2015). However, it is far more common for researchers to study broader student attempts at proof construction and to categorize the features of these attempts in the analysis phase of the research (e.g., Bell, 1976; Weber 2001). If patterns emerge in these attempts, conclusions and frameworks can then be crafted to reflect commonalities in successful attempts, failed attempts, and persistent student strategies.

Studies that highlight failed proof attempts emphasize common areas that prevent students from arriving at a completed proof (e.g., Weber, 2001; Moore, 1994). Of note, Moore (1994) identified seven categories of such student difficulty by examining previous proof literature: (1) failure to state the definitions; (2) failure to understand the meaning of the concepts intuitively; (3) failure to use the concept images while performing proof; (4) lack of generalization and use of examples; (5) failure to know what kind of proof structure to use from the definitions; (6) failure to understand the mathematical language and notations; and, (7) failure to know how to begin the proof process. Edwards & Ward (2004) expanded on Moore's identification of definitions as source of difficulty for students by focusing on students' use of definitions in their proof attempts in addition to exploring student beliefs. Although the definitions needed to approach the tasks were available to participants during the course of the task, Edwards & Ward still identified four issues prevalent among students' task responses that interfered with proof construction: (1) Many students do not *categorize* mathematical definitions the way mathematicians do; (2) Many students do not *use* definitions the way mathematicians do, *even when students can correctly state and explain the definition*; and (3) Many students do not use definitions the way mathematicians do *even in the apparent absence of any other course of action* (pp. 415-418).

Despite the persistent identification of proof difficulties in research with undergraduate mathematics students, research demonstrates that even if these difficulties are supposedly resolved by student possession of conceptual knowledge, this does not necessarily result in a student being able to write proofs (Weber, 2001). Weber examined graduate and undergraduate students that he deemed as possessing the conceptual knowledge and ability to apply facts needed to prove a statement. However, the graduate students in his sample were able to construct normed proofs, while the undergraduate students were typically unable to prove the statements provided. Weber contrasted the graduate students' success with the failure of undergraduate students to construct a proof and concluded that the key difference between the graduate students and the undergraduate students was that doctoral students knew which theorems and techniques would be useful, as well as when symbolic manipulation would be the most efficient way to tackle a problem. Weber refers to this ability as strategic knowledge.

Zazkis, Weber, & Mejia-Ramos (2015) deliberately studied students they deemed to be "highly successful mathematics majors" to determine the strategies common to these students, and if these strategies could also ensure success for mathematics majors who adhered to them. By presenting novel tasks that did not require more than familiarity with first semester Calculus concepts to these students, the authors found that these students were prone to engage in what the researchers described as the "targeted" strategy and the "shotgun" strategy when creating proofs. In the targeted strategy, students "would develop a strong understanding of the statement they were proving, choose a plan based on this understanding, develop a graphical argument for why the statement is true, and formalize this graphical argument into a proof" (Zazkis, Weber, & Ramos, 2015, p. 15). When using the shotgun strategy, students would "begin trying different proof plans immediately after reading the statement and would abandon a plan at the first sign of difficulty" (Zazkis, Weber, & Ramos, 2015, p. 16).

In the shotgun strategy, students' rapid dismissal of certain techniques was an attempt to find the quickest way to solve the problem. If too much work was involved, students supposed that there must have been a better way to approach the problem (this could possibly be identified as a result of the individual student's beliefs.) The targeted strategy, on the other hand, is described as requiring in depth understanding of the task at hand to determine why the required result holds true. Although the authors do not claim that these techniques are requirements to be successful in higher level mathematics, the six students studied that "successfully" completed proof tasks overlapped in employing these strategies.

The shotgun strategy as previously described is analogous to the brute force method as described by Duffin & Simpson (2006) in which the authors reasoned that these sort of brute force methods require a substantial amount of exposure to various techniques, which leads to the implication that students must see enough proof strategies in order to be successful. Duffin & Simpson found that this technique becomes ineffective, or results in failure to prove, as students advance in their mathematics education, ultimately being abandoned entirely once students enter graduate mathematics courses (Duffin & Simpson, 2006). The shotgun, or brute force method, is at odds with what is required of students in advanced mathematics programs.

### *Student and instructor beliefs about proof.*

Proof-specific beliefs research focuses on participants' beliefs about proof and can offer insight into why participants make specific choices in the act of proving (e.g., Knuth, 2002; Weber & Mejia-Ramos, 2014; Balacheff, 1988; Brown, 2017; Stylianou, Blanton, & Rotou, 2015). This type of research typically appears in three forms: (1) student beliefs about what constitutes a proof; (2) beliefs about the purpose of proof (e.g. Healy & Hoyles, 2000); and (3) beliefs about what makes a proof convincing (Knuth, 2002). Although there are many studies addressing these kinds of beliefs, the literature relevant to this study compare the beliefs of collegiate instructors and/or mathematics education researchers to students, prospective teachers, and practicing teachers.

At the collegiate level, mathematics professors have indicated that they present proofs to mathematical statements in their mathematics courses under the assumption that students can follow and understand the proofs, but, when these mathematics majors are asked to construct a proof of the same mathematical statements on their own, they resort to restating the words used by their professor when reciting the proof (Weber & Mejia-Ramos, 2014). Weber and Mejia-Ramos also found that many students hold some combination of four beliefs about proof reading and the purpose of proof that are "unproductive": (1) students believe that they should not be expected to justify statements within a proof that has been given to them. If a statement is given within a proof, then students felt that it should be explicitly stated how that assertion follows from previous assertions; (2) students believe that understanding a proof consists entirely of understanding how each statement in the proof is justified; (3) students do not think that studying a proof is a lengthy process; and (4) students believe that they should not be expected to draw diagrams to help them understand a proof. If a diagram would aid a students' comprehension of a proof, then it should be provided within that proof.

These beliefs were in opposition to what professors believed was necessary for students to understand proofs (e.g. proof constructions take time, drawing pictures are helpful, etc.) The results from the Weber and Mejia-Ramos' study suggest that this conflict of beliefs results in mathematics majors not meeting the self-sufficiency and time expectations of their professors, thus leading to a failure to attain the mathematical goals expected in mathematics courses. Student mimicking professor explanations may indicate that student over-reliance on recall may be a source of student difficulty with proof construction.

Knuth (2002) utilized five common perspectives among mathematics educators on the role of proof within mathematics to consider secondary mathematics teachers' conceptions of proof. The five roles that proof is described as playing in mathematics are: (1) verification of a statement as true; (2) explication of why a statement is true; (3) mathematical communication of knowledge; (4) to create new mathematics through discovery; and (5) to form axiomatic systems through systematizing statements (Knuth, 2002). Knuth found that beliefs about the role of proof among teachers were similar to the five common perspectives. The participants believed that the roles of proof are: (1) establishment of truth; (2) explanation; (3) communication of mathematics; and (4) creation of knowledge/systemization of results, however, these beliefs were specific to the secondary mathematics context (Knuth, 2002, p. 386). None of the teachers in this study indicated that the role of proof was to develop understanding, and teachers relied on the truth of a statement to determine if the argument they were presented was a proof. When uncertain about the truth of a statement, participants used empirical evidence to establish conviction.

Of particular importance to my dissertation study is the examination of teacher's methods for evaluating arguments. The teachers in this study were given five propositions and given between three and five proof constructions justifying each of these propositions (Knuth, 2002, p. 391). The methods utilized in proof constructions were varied. Participants were asked to rate the validity of each argument on a scale of one to four[4] where one represented an argument that teachers deemed not-a-proof and invalid, ratings of two and three were included to allow teachers to express varying levels of belief that the argument was a proof or not, and a rating of four represented the belief that the argument was valid. While participants displayed great strengths in identifying valid arguments, they also rated many non-proofs as proofs (Knuth, 2002, p. 391).

The criteria used by participants in evaluating the arguments were: (1) valid methods as evidenced by the type of "formal" method used rather than the reasoning employed; (2) mathematically soundness as evidenced by the reasoning within the proof without regard for the method; (3) sufficient detail as evidenced by the participant's ability to follow the argument; and (4) knowledge dependence as evidenced by the specialized knowledge needed to understand the argument. Detail and knowledge dependence were used to classify arguments as proofs or not, while the perturbations in ratings were based on the teacher's perception of method validity and mathematical soundness. As a final comment, the most common characteristics that made a proof convincing to participants were: (1) concrete features (i.e., using specific values of a visual representation; (2) familiarity; (3) sufficient level of detail; (4) generality; (5) the proof "shows why"; (6) valid method; and that (7) all arguments in the set were equally convincing (Knuth, 2002, p. 398). I anticipate that many of the participants' criteria for the evaluating proof constructions and proof validity will be in line with the methods of the student participants in my dissertation study.

### *Theories of proof.*

Harel & Sowder's (1998; 2007) theory of proof schemes (e.g. Lee, 2015) provide a valuable framework for determining where a student may be in the process of individualization, though this process need not be linear and a student may possess several schemes at once. Specifically, the authors state that these schemes represent a stage of cognitive development that corresponds to intellectual ability. A person's proof scheme consists of what constitutes ascertaining and persuading for that person (Harel & Sowder, 1998; 2007). The diagram created by Harel & Sowder (2007) that represents their identification and categorization of types of proof schemes can be found in Figure 3.

---

[4] Knuth's (2002) choice of a four-point scale provides precedent support for my use of a four-point scale with instructor participants in my study. In contrast, I allow the instructor participant to define their four-point scale according to their own narrative preference and I do not make determinations about the validity of student proof constructions. Further discussion of narrative ratings can be found in my methods.

**REDACTED IMAGE**

Figure 3: Diagram of proof schemes. Reprinted from Harel, G., & Sowder, L. (1998). Students' proof schemes: Results from exploratory studies. In A. Schoenfeld, J. Kaput, & E. Dubinsky (Series Eds.), *Issues in Mathematics Education: Vol. 7. Research in collegiate mathematics education. III* (pp. 234-283). doi: 10.1090/cbmath/007

The three main categorizations of proof schemes are the external, empirical, and analytical proof schemes. *External proof schemes* are defined as those that represent a conviction outside of the student and they fall into three categories; *authoritarian*, *ritual*, and *symbolic* proof schemes. Of these schemes, Harel & Sowder state, "When formality in mathematics is emphasized prematurely, students come to believe that ritual and form constitute mathematical justification. When students merely follow formulas to solve problems, they learn that memorization of prescriptions, rather than creativity and discovery, guarantee success" (Harel & Sowder, 1998, p. 246).

The external proof schemes can be summarized (based on what is used to convince oneself) as follows:

- Ritual – the ritual of argument presentation. An acceptable proof is one that looks like what a student expects of proof
- Authoritarian – the word of an authority (lecturer, textbook, etc). Mathematics is seen as finished truth on which students can have little impact.
- Symbolic – the symbolic form of the argument. Symbols are manipulated without reference for meaning or parameters regarding how they should be used

The second category of proof schemes are the *empirical*. These proof schemes are "validated, impugned, or subverted by appeals to physical facts or sensory experiences" (Harel & Sowder, 1998, p. 252) and are split into two categories; inductive and perceptual. Students with inductive proof schemes convince themselves through quantitative verification; these students will make mathematical arguments that check specific arguments in one or more instances. If a pattern appears to these students, their belief in the truth of the argument will strengthen. Of perceptual proof schemes, the authors state, "Perceptual observations are made by means of rudimentary mental images- images that consist of perceptions and a coordination of perceptions, but lack the ability to transform or to anticipate the results of a transformation. The important characteristic of rudimentary mental images is that they *ignore transformations on objects or are incapable of anticipating results of transformations completely or accurately*" (Harel & Sowder, 1998).

The final and largest proof scheme category is the *analytical* proof scheme. This can be described as a proof that "validates conjectures by means of logical deductions" (Harel & Sowder, 1998). The two main subcategories of analytical proof schemes are transformational proof schemes and axiomatic proof schemes. While these proofs schemes have sub-categories, for the purposes of this literature review I will not break down these subcategories further. Transformational proof schemes are characterized by student ability to anticipate results. Student behavior in this scheme has *intent* and "involve operations and anticipations of the operations' results" (Harel & Sowder, 1998, p. 258) and this is the proof scheme that students are expected to move toward as they mature in mathematics. Additional hallmarks of this proof scheme are the ability of students to mentally transform images and consideration of the generality of their arguments. Axiomatic proof schemes are those "when a person understands that at least in principle a mathematical justification must have started originally from undefined terms and axioms (facts, or statements accepted without proof)" (Harel & Sowder, 1998, p. 273). The limitation of this proof scheme as described by the authors is, "He or she, however, may be able to handle only axioms that correspond to her or his intuition, or ideas of self-evidence…" (Harel & Sowder, 1998, p. 273).

Most students are thought to possess an authoritarian proof scheme, and Harel & Sowder (1998) cite schooling as influential in students developing such schemes. They state, "This begins with elementary mathematics where children are rushed into using mathematical prescriptions to solve arithmetic problems (Harel, 1995) and continues with secondary and postsecondary mathematics where instrumental understanding rather than relational understanding is emphasized throughout the curriculum" (Harel & Sowder, 1998, p. 247). This student approach to mathematics can therefore fall into what Schoenfeld (1991) refers to as non-sensemaking; students performing mathematics with little regard for meaning.

Using Harel & Sowder's (2007) proof schemes framework, a student engaging in a targeted strategy (Zazkis, Weber, & Mejia-Ramos, 2015) would likely be categorized as possessing an analytical proof scheme, while students using brute force and shotgun methods (Duffin & Simpson, 2005) would likely be categorized as possessing authoritative proof schemes. However, it is difficult to ascertain if either of these students' discourse around the *concept* assessed is normed. This is important because Harel & Sowder's theory is explicit in its identification of expert discourse of proof as ideal. Students categorized as having transformational proof schemes in the study were only assigned this categorization based on the researcher's interpretation of student actions. The authors do not explore the students' sense-making while they attempt the tasks, but rather, they attempt to categorize student choices in their proof construction. Theorized shotgun and targeted strategies require student exposure to an extensive amount of proof examples, and this could plausibly lead to a student proof scheme that is more ritualistic (by Harel & Sowder's definition) in nature than analytical. In short, while Weber demonstrated that students with content "knowledge" may still struggle with proof construction, I argue that it is also possible that students are able to engage in successful proof construction without content knowledge. The role of non-proof content knowledge and how it appears within each scheme are not explored within the Harel & Sowder framework.

What is interesting in the proof schemes as defined by Harel & Sowder is the separation of the external proof schemes from that of the analytical. Although Harel & Sowder state that it is possible for students to hold multiple proof schemes at the same time, Sfard's theory of commognition can perhaps be used to illustrate the ways that students who actually hold external proof schemes can be incorrectly categorized as possessing analytical proof schemes, particularly if a student has robust experience and exposure to proof that would help them finely tune their personal set of *applicability conditions*; defined as the rules that "delineate, usually in a nondeterministic way, the circumstances in which the routine course of action is likely to be evoked by the person" (Sfard, 2008, pg. 209) using the student's own set of clues within key words from a specified *task situation.* Sfard's idea of routine, particularly Sfard's definition of ritual and deeds, may still be present in students categorized in each of these proof schemes. In particular, analytical proof schemes may still be ritualized in nature if students can successfully coordinate multiple rituals to construct proofs that meet the expectations of their course instructors.

### Research on injective and surjective functions.

A focus on success and failure with proof construction can obscure the extent of actual student engagement with content. Injective and surjective maps are prevalent within upper-level mathematics courses like Abstract Algebra and Real Analysis, and students must lean on the discourses around function that they enter these courses with to continue successfully in mathematics. However, very few studies have been conducted on student "understanding" of injections and surjections (Thoma & Nardi, 2020; Bansilal, Brijlall, & Trigueros, 2017), possibly because students are able to engage in ritualized performance on such tasks to the satisfaction of their professors. Far more attention has been paid to the parent topic of functions, and much of the information available about student difficulty with injections and surjections is a by-product of research on student "understanding" of functions (e.g., Moore, 1994; Breidenbach, et.al., 1992). Of note, functions present a profound difficulty for high school students and, in many cases, even their teachers are not fully aware of the norms within function discourse (Briendenbach, Dubinsky, Hawks & Nichols, 1992; Dubinsky & Wilson, 2013). Naturally, it follows that a student who struggles with function discourse will most likely struggle with discourse on injective and surjective functions.

As compiled by Selden & Selden (1992), students exhibit considerable difficulty within function discourse due to the following common student object-level rules that are in conflict with mathematics community norms: (1) functions must be a formula; (2) a change in independent variable is required; (3) functions should be defined by a single rule; (4) non-algebraic functions must be approved by the mathematical community; and (5) a function must be one-to-one. The belief that the identification of a relation as *one-to-one* has bearing on if that relation can be called a function can also appear among students because of the common student metarule that the definition of one-to-one is what constitutes the definition of a function. In short, students tend to not only believe that all functions can be named one-to-one, but that being named one-to-one is what allows a relation to be named a function. Furthermore, Selden & Selden explain that a student may fail to reconcile the mathematical definition of a function, the definition of a function as recalled by the student, and the student's personal metarules for identifying a function (Selden & Selden, 1992).

Furthermore, the realizations of a student can be heavily dependent on how a function is signified. The common representations of a function are: (1) ordered pairs; (2) explicit equation; (3) graphical; and (4) regions mapped by arrows (Selden & Selden, 1992). Tirosh & Tsamir (1996) found that when students were tasked with determining the size of sets, functions represented in numerical vertical form often resulted in students appealing to one-one correspondence to conclude that sets were of the same size. When presented with alternate representations, students often struggled to arrive at a conclusion about set size. This phenomenon is also present when students are tasked with determining whether they have been presented with a function. Scheja & Pettersson (2010) argue that for a student to advance their discourse on functions, the ability to switch between contexts, and therefore representations, plays a critical role.

Dubinsky & Wilson (2013) describe a student as possessing "understanding" of function if they can recognize functions and one-to-one functions given multiple representations. If Dubinsky & Wilson's criterion for "understanding" functions is acceptable, this suggests that many students, through lack of introduction in their mathematics classrooms, are limited in their ability to realize varying signifiers of relations as functions. Norman (1993) found that the secondary education teachers in his study consistently appealed to graphical representations of functions in response to problem tasks, and that these teachers almost exclusively used the vertical line test to determine if the image in the graph represented a function. It should be of no surprise then that when students are asked to supply their definitions of a function, many claim that passing the vertical line test suffices, neglecting to consider domain and range (Wilson, 1994). Other features prominent in student-contributed definitions of a function are (1) one-valuedness; (2) lack of discontinuity; (3) lack of a split domain; and (4) lack of exceptional points (Vinner & Dreyfus, 1989).

For a student with any of the previously described discourses to build and utilize a normed discourse of onto and one-to-one within function, they must encounter task situations that force them to either abandon or add to their discourse on functions due to their discourse being at odds with the one currently expected of them (Sfard, 2008). Sfard describes any situation, whether interpersonal or intrapersonal, where communication is attempted across incommensurable discourses as *commognitive conflict*. Student or instructor awareness of this conflict should prompt an attempt to resolve the conflict. If a student whose discourse on function is not normed is never confronted with a task situation that introduces conflict, there is no impetus for a student to alter their discourse. I will now consider two definitions of one-to-one and onto.

Smith, Eggen, & St. Andre (2014, pp. 205-208) defines one-to-one and onto functions in the

---

(1a) A function $f: A \to B$ is onto B (or is a surjection) iff $\text{Rng}(f) = B$

(2a) A function $f: A \to B$ is one-to-one (or is an injection) iff whenever $f(x) = f(y)$, then $x = y$

---

following manner:

while Chartand, Polimeni, & Zhang (2017, pp. 256-257) provide the following definitions:

> (1b) A function $f: A \to B$ is **onto** or **surjective** if every element of the codomain $B$ is the image of some element of $A$
>
> (2b) A function $f: A \to B$ is **one-to-one** if whenever $f(x) = f(y)$, where $x, y \in A$, then $x = y$

In the terminology of commognition, these statements are expert endorsed narratives about injective and surjective functions. A function $f$ is *called* a surjective function from A onto B such that when $f$ is evaluated at every point in A, the resulting collection of objects forming the range *is* B. This can alternatively be described by stating that for every y in B, there exists some x in A such that $f(x) = y$, and this relationship need not be unique. To substantiate the identification of a function as injective, students must recall that if the values in the range (Rng) B of a function at two points of the domain A are the same, then the points selected from the domain must be the same point. To paraphrase, whenever f is evaluated at an arbitrary point in its domain, it maps to a unique point in the range. The domain and codomain are not restricted to the field of the naturals, integers, the reals, etc., but a collection of non-numerical objects can also serve as the domain and the codomain. For example, students at a university may be assigned a unique identifier consisting of eight digits, where the collection of the names of all students are in the domain and all assigned identifiers are in the range. There is a process by which, for each of these students' names, a unique eight-digit identifier is generated. If two supposed students appear to have the same identifier, then these students must be the same student.

Moore (1994) provides the following equivalent definitions of injective (one-to-one) functions:

(1) A function f is one-to-one if and only if for all x and for all y in the domain of f, if f(x) = f(y), then x = y

(2) A function f is said to be a one-to-one function provided that no two distinct members of f have the same second term

The first definition is the most prevalent within collegiate classrooms and textbooks and is the one provided by Smith, Eggen, & St. Andre (2014) and Chartand, Polimeni, & Zhang (2017). According to Moore, this definition provides students with a technique for proving that a function is one-to-one. This technique involves evaluating f at two arbitrary points and then deriving an equality, or disproving that a function is one-to-one by deriving a contradiction. Moore states that the second definition (the ordered pair perspective) may be difficult for a student to reconcile with the first. I note that both textbooks authors and Moore provide definitions for one-one that differ slightly in their appearance when compared to the other sources.

Moore asserts that students sometimes "know" a definition but cannot use that definition within a proof. He provides an example of a student who was able to use the second definition to provide examples and nonexamples of one-one functions, but was who unable to prove that a function should be called one-one, even when provided with the first definition. Additionally, Moore found that students experienced difficulty because they did not understand the hypothesis of a specified one-to-one problem, pointing to yet another of Moore's identified struggles with proof that is unrelated to content itself: inability to understand and use mathematical language. From a commognitive perspective, students who struggle in these ways likely do not associate any realizations or routines to the signifiers present in the definition.

    These prior studies in the realm of functions, injections, and surjections illustrate that how a problem is posed, and potentially the signifiers and realizations within or not within the task situation that can serve as key words, can be just as much of a hurdle for students in proof construction than lack of acquaintance with definitions. Even with such an acquaintance, students may still struggle with proof construction despite the explicit availability of an exploitable definition. These oddities were useful in informing the crafting of the methods necessary for my study.

## Chapter 3: Methods

Recall that I seek to answer the following questions:

1. In what ways do undergraduate mathematics students interpret and make sense of injection, surjection, and proof instruction in transition to proof courses?
2. In what ways do instructors interpret and make sense of proof, injective functions, and surjective functions?
3. How do the injection and surjection discourses of students differ from the discourses of their instructors?

The commognitive researcher is unique among participationists in viewing mathematics as discourse. Positioning discourse as the unit of analysis allows the commognitive researcher to minimize the amount of interpretation, objectification, and subsequent subjectification of participants in conducting, analyzing, and reporting on research. Although the intent of this positioning is to produce research where the integrity of participants' utterances is preserved and the voice of the researcher is minimal and obvious when used, commognitive researchers acknowledge that the very presence of a researcher can impact utterances and activity. With this limitation in mind, Sfard (2008) asserts that commognitive research must follow five methodological principles:

- Principle of operationality – the researcher must attend to clarity and protect themselves against misunderstandings. Communicational effectiveness is accomplished though the operationalization of any terms with objectified uses. For example, knowing, learning, and understanding are terms with highly varied uses within mathematics education research. They each describe unknowable states, and therefore they must be operationalized before their use in commognitive research. Any frameworks or codes that I create must be operationalized.
- Principle of completeness – The target of study should be the entire discourse surrounding the target. This study uses injective and surjective functions to explore sense-making in proof construction. Therefore, mathematical discourse is the unit of analysis, though I ask participants to represent a portion of that discourse. Student

- utterances during interviews are relevant to the analysis without regard for if these utterances refer to proof, injective, or surjective functions specifically.
- Principle of contextuality – All manner of interaction is a learning event. Data must be presented within its context. The commognitive researcher must attempt to "document human interaction fully" (Sfard, 2012, p. 8).
- Principle of alternating perspectives – The commognitive researcher must act as both an insider and an outsider when analyzing data. It is possible that the researcher and the participant have incommensurable discourses. Stepping into the role of outsider allows the researcher to eliminate interpretation that is specific to their own experiences and biases. Anyone who consumes the research can then make their own conclusions.
- Principle of directness – The commognitive researcher should report their results in a direct fashion without editorializing. The commognitive researcher begins by reporting was what said or done by the participant first. The narratives presented by the researcher should be identified as not being "about the world, but about the participants' narrative about the world" (Sfard, 2012, p. 8). Any inferential narratives that a researcher makes in turn are also not about the world, but about the researcher's narrative about the participant's narrative about the world.

These principles are consistent with my proposed methods for answering my research questions. The results of an initial feasibility study informed the methods to be described herein and details concerning the feasibility study can be found in Appendix H. I will use a combination of field notes, semi-structured clinical interviews (Clement, 2000; Patton, 2015), task-based interviews (Goldin, 1997, 2000; Maher & Sigley, 2014), and an adaptation of concept maps (Novak & Cañas, 2008; Novak & Gowin, 1984; Novak, Gowin & Johansen, 1983; Williams, 1998) over the course of this study. The primary goal of Phase 1 of my study is to produce artifacts that can aid me in anticipating Phase 2 participants' evoked discourse and build models of typical student solution paths. The goal of the initial clinical interview with Phase 2 participants is to produce

| Study Phase | Intended Timeline | Alternative Timeline (Phase 2 occurs in Fall 2022) |
|---|---|---|
| Phase 1 – Individual student phase and ongoing analysis | Late February/Early March – April, 2022 | Late February/Early March – April, 2022 |
| Phase 2 – Class-based phase with instructor and their consenting students | Late May – early August 2022 | Late July – early December 2022 |
| Phase 3 – Final clinical interview with instructor and whole study data analysis | Late July – August 2022 | December 2022 |
| Phase 4 (Optional) – Repeat Phases 2 and 3 | Early August – December 2022 | January – May 2023 |

data that can approximate the course instructor(s) and student participants' function discourses prior to any instruction in the introductory proof course. I will then use task-based interviews throughout the semester for the remainder of this phase. During Phase 3, which takes place after the conclusion of the semester, I will have a final clinical interview with the course instructor. In this chapter, I will provide my rationale for focusing on students in transition-to-proof courses specifically, the university setting of the study, the methods of data collection, and my proposed data analysis subject to the five principles of conducting commognitive research. I include two possible timelines for this study contingent on when the class-based phase will be conducted.

**Participants and Setting**

The participants in this study are: (1) undergraduate students enrolled in an introductory proof course intended for mathematics majors without regard for the instructor of the course, (2) the instructor(s) for an introductory proof course, and (3) the instructor participant's consenting students. Phase 1 of the study is not attached to a classroom nor instructor and I will recruit participants by advertising the study in UGA Aderhold Hall, by emailing students through a departmental e-mail list upon receiving permission of the department, and/or by requesting that an instructor of a transition to proof course send the study advertisement to their students. Students that participate in this phase of the study must be enrolled in a transition to proof course and 3-7 students are expected to participate in this portion of the study. Phases 2 and 3 of the study focus on an instructor(s) and the students in their course. I will recruit participants by advertising the study in an email via the course instructor. Additionally, I will ask instructors to advertise the study in their classes. The selection of the instructor(s) for the remaining phases of the study is dependent on the typicality of course instruction as evidenced by past course syllabi when compared to syllabi written by instructors at multiple institutions. Past introductory proof syllabi should clearly delineate the intended course schedule and content, course format, instructor goals for the course, expectations for students, and the intended course textbook(s). Typicality will be assessed based on the stated content to be taught and course format among sample syllabi.

In the event of limited instructor options for the class-based portion of the study, I will conduct an additional class-based phase (Phase 4) to meet the goal of typicality of course instruction. This will strengthen the overall aims of my dissertation by yielding comparative data for analysis between students receiving typical instruction and students receiving atypical course instruction. An example of instruction that is currently atypical in transition to proof courses would be a course delivered in a flipped or active learning format as opposed to the traditional lecture format. Single and two-semester course instructors will still be considered typical if all other conditions for typicality are met. All students in the selected instructor's introductory proof course will be invited to participate in the initial training and discursive map portion of the class-based study. Students will be invited to continue in the study upon verification that the student can construct a discursive map of fraction according to the training guidelines. I would like to invite ten students to the initial clinical interview and at least six to the second interview. Ideally, three students would complete the fourth interview.

**Discursive Maps**

**Concept maps.**

Given the limited historical literature on using realization trees to study discourse development (e.g., Sfard, 2008; Caspi & Sfard, 2012; Newton, 2012; Heyd-Metzuyanim, 2012;

Weingarden, Heyd-Metzuyanim & Nachlieli, 2019), I sought alternative support for their use in my study within the mathematics education literature. Sfard's (2008) perceptually accessible and hierarchal commognitive *realization trees,* as well as their symbolic artifact encompassing *symbolic tree* counterpart, are similar in appearance and features to concept maps (Novak & Cañas, 2008; Novak & Gowin, 1984). Concept maps were developed as a way to visualize knowledge structure to assist both students and instructors in developing their "understanding" of a concept (Novak & Gowin, 1984). Recall that the term *concept* has been defined discursively as a word or signifier together with its uses (Sfard, 2008) within the commognitive framework.

 Concept maps are created by placing words into a diagram and indicating connections between these words by drawing a line between them. These lines are typically labeled using linking words (Novak & Gowin, 1984). When two objects are read together with their linking word, these should form meaningful propositions to the author of the map, and the overall map should be hierarchical. These maps were developed within a cognitivist view, though Novak & Gowin made significant reference to the power of language in describing regularities and making communication effective in their description of concept maps. The linguistic focus of Novak & Gowin strengthened concept maps' compatibility with the commognitive lens for analyzing changes in student discourse and differences between student and instructor discourse.

 Concept maps have been used to study differences in function "concept" between instructors, traditional students, and reform students in a Calculus course (Williams, 1998), to study differences in function knowledge structure between high and low gain students

**REDACTED IMAGE**

Figure 4: Algebra concept map constructed for math review course. Reprinted from Novak, J. D., & Gowin, D. B. (1984). *Learning how to learn,* p. 78. New York, NY: Cambridge University Press.

(McGowen & Davis, 2019), and to assess engineering competence (Watson, Pelkey, Noyes, & Rodgers, 2016). When examined alongside traditional measures of student performance (e.g. SAT), Novak, Gowin & Johansen (1983) found that concept maps were better able to account for the extent of student engagement with concepts given students' autonomy in what they choose to include and represent about the concept in their maps. However, they also found that concept maps and course grades showed no correlation, while SAT scores and course grades were highly correlated. This could possibly indicate that SAT scores and course grades are measures of similar features of student activity; based on the traditional features of the SAT, I hypothesize that these are features that can be assessed as "correct" vs. "incorrect" or "complete" vs. "incomplete". As a reminder, these maps were required to be hierarchal. Of note, numerous studies on concept maps established their validity for documenting and exploring conceptual change (Williams, 1998, p. 417). Within the domain of mathematics specifically, Park (1993) found a strong correlation between concept-map scores and post-test scores among students in a collegiate calculus course. As a note, concept-map scoring (Watson, Pelkey, Noyes, & Rodgers, 2016) methods can widely vary and will therefore impact correlations with traditional measures of achievement. For the purposes of this study, I do not wish to create a scoring system, which is typically reliant on a hierarchy. I am interested in assessing concept map coherence relative to the participant's own discourse.

### Adapting concept maps for commognitive analysis.

      For the purpose of adapting concept maps for use in this dissertation, two key features of concept maps required modification. First, Novak's envisioned concept maps are hierarchical. Within mathematics, however, numerous mathematical symbols, words, and representations are regularly interchangeable under specific conditions; though a different set of conditions can either restrict or expand the set of interchangeable mathematical objects. This makes the notion of hierarchy within a concept map problematic for its use within mathematics. Furthermore, as concepts are developing (e.g. as in the case of students within a course) it should not be expected that the author can identify a hierarchal relationship between words. Even if I were to dispense of hierarchy within concept maps due to the precedent and varied implementations of concept maps within mathematics education research (e.g. Williams, 1998; McGowen & Davis, 2019), I would also need to justify removing the requirement for hierarchy within realization trees. To find support for or against this decision, I revisited Sfard's original theory of commognition.

      Though Sfard (2008) and I agree on the point that symbolic artifacts are regularly interchangeable, and therefore, a signifier-realization dichotomy often fails to exist (p. 164), Sfard does not make any reference to a feature or lack thereof regarding hierarchy within symbolic trees as opposed to realization trees. After the symbolic tree terminology is introduced in *Thinking as commognition* (Sfard, 2008), Sfard refers to realization trees and symbolic trees interchangeably throughout the remainder of the text; this would seemingly indicate that symbolic trees should therefore follow the same properties of realization trees. Additionally, Sfard remarks that realization trees rely on an evoking situation and are therefore partial; she categorically states that a complete realization tree for a participant would be untenable for a researcher to create. Without evidence to the contrary, I expect this sentiment to extend to symbolic trees as well.

      I re-directed my focus to locating implicit support from Sfard for removing the hierarchal requirement for realization trees. Sfard specifies that, "…whether a word, algebraic symbol, or

icon should count as a signifier or as a realization of a signifier is a matter of use, not of any intrinsic property of these artifacts" (Sfard, 2008, p. 164). However, seemingly in support of the manner that I want to use realization or symbolic trees to study discourse, Sfard mentions that one way to explore a person's discourse is by examining all of a person's uses for a signifier. This is, in fact, Sfard's definition of a mathematical object. Recall that Sfard states, "*mathematical objects are abstract discursive objects with distinctly mathematical signifiers*" (Sfard, 2008, p. 172 emphasis in original) and that the discursive object signified by S in a

**REDACTED IMAGE**

Figure 5: Compound phrase signifier and a possible partial realization tree. The phrase is unpacked through a series of skillful realization procedures to determine the object of talk. Reprinted from Sfard, A. (2008). *Thinking as communicating: human development, the growth of discourses, and mathematizing.* Cambridge University Press.

discourse on S is defined as, "the realization tree of S within this discourse" (Sfard, 2008, p. 166). As a final piece in the justification of removing the requirement for hierarchy in realization trees, Sfard notes the following about realization trees:

> "First, realization trees, and hence mathematical objects, are *personal constructs*, even though they originate in public discourses that support only certain versions of such trees. As researchers, we may try to map personal realization trees and present them in diagrams…Second, the realization trees are a source of valuable information about the given person's discourse. Making skillful transitions from one realization to another is the gist of mathematical problem solving. In addition, a person's tendency to apply mathematical discourse in solving practical problems depends on her ability to decompose signifiers into trees of realizations with branches long enough to reach beyond the discourse, to familiar real-life objects and experiences…This last statement leads me to a third point. While analyzing transcripts of conversations in an attempt to map discursive objects, one needs to remember that these *personal constructs may be highly situated* and, in particular, can be easily influenced by interlocuters and by other specifics of the given interaction." (Sfard, 2008, p. 166, emphasis in original)

Sfard's original expectation of hierarchy is borne of Sfard's assertion of mathematics as hierarchal; within normed mathematical discourse I am inclined to agree. If I create the realization tree on behalf of the participant, I can use my role as an insider to expert mathematical discourse to structure the tree in a hierarchal manner, however, that structure may

**REDACTED IMAGE**

Figure 6: Partial realization tree for *basic quadratic function*. Reprinted from Nachlieli, T., & Tabach, M. (2012). Growing mathematical objects in the classroom -The case of function. *International Journal of Educational Research*, *51-52*, 10-27.

not necessarily reflect how the participant uses the signifier. If a participant does not depict or indicate their uses of a signifier in a hierarchal manner, then I have no reason to assume that specific discourse is hierarchal. Inferring otherwise would violate the principles of commognition.

As a result of my interpretation of Sfard's utterances, it appeared fitting to: (1) retain the hierarchal structure of realization trees in depictions of student problem solving as this represents the participant's successive parsing of the problem. I note that Sfard intended for realization trees to be recursively unpacked from a root node where each successive realization layer retained all of the narratives about the parent signifier, so I refer to structures that are faithful to Sfard's realization and symbolic trees as *evoked trees* or *trees* within this study; (2) to approximate a participant's baseline discourse for a signifier within their mathematical discourse (i.e. the participant's uses of the mathematical signifier), I shift the task of creating a symbolic/realization tree to the participant (Sfard, 2009) and remove the requirement of hierarchy. Of course, participants can still choose to represent their discourse in a hierarchy, but their tree will still be acceptable without it; and (3) to minimize the situatedness of the baseline discourse created for a signifier, the participant is provided with only the signifier as a prompt. Trees created in the absence of a specific task are referred to as *unrestricted* in line with Williams' (1998) use of the term unrestricted when referring to concept maps.

Interestingly, a study of unrestricted function concept maps within Williams (1998) demonstrates that experts were inclined to hierarchal structuring of mathematical content without prompting, while hierarchy was largely absent among student concept maps. With the removal of the requirement of hierarchy due to developing discourses among students, the student's agency in creating their own maps as opposed to the researcher, and the lack of a need for hierarchal structure for scoring purposes, I considered the first issue with adapting concept maps resolved.

**REDACTED IMAGE**

Figure 7: Partial hierarchal realization tree for function. I reimagine this map as an unrestricted discursive map. Reprinted from Newton, J. (2012). Investigating the mathematical equivalence of written and enacted middle school standards-based curricula: Focus on rational numbers. *International Journal of Educational Research*, *51-52*, 66-85.

The second issue that needs to be resolved within concept maps is the lack of a discursive context. It is therefore necessary to appeal to a discursive definition of concept to closely align concept maps to the commognitive framework. Therefore, concepts within the map are read as signifiers together with their uses. In short, concepts must be considered as mathematical objects consisting of a word or symbol together with its uses. I retain the feature of allowing linking words within concept maps due to allowing the participants in the study to author their own map. Linking words provide further clarification on the uses of a signifier that the participant intends to represent. With these two modifications (optional hierarchy and a discursive definition of concept), I differentiate the maps I will use in this study from the original concept maps as well as realization trees by referring to the unrestricted, discursive form as *discursive maps* or discursive concept maps (where concept is defined discursively or on a purely linguistic basis).

While concept maps are not seen as "windows into the mind of the author" (Novak & Gowin, 1984), I intend to use the discursive map version in concert with student created and hierarchical evoked trees that are specific to the tasks-situations they encounter to better approximate student thinking. Due to the construction of the trees being situated in the task-situation, I can construct rich descriptions of how an individual student's discourses are organized along with how portions of these discourses are activated and accessed in specific situations. To reiterate, this latter feature of the study is an acknowledgement of the impact of situatedness on students' mathematical thinking and activity. In summary, discursive maps are an unrestricted, hierarchy optional, scenario-independent form of realization trees where the root signifier is a mathematical object together with its uses. Participants in this study are asked to discuss the function object only. The evoked trees that students will use to describe their successive series of realizations are the students' attempts to communicate their interpretation of a task situation. They are not asked to discuss all of the uses of a mathematical object in response to the task-situation. Rather, they indicate the specific uses that they have considered in their reading of the task.

**Discursive map validity.**

After I established support for the removal of hierarchy from concept maps and discursive maps, I initially focused on establishing that concept maps could be parsed into commognitive items. Focusing on two studies, I collected ten function concept maps from seven unique participants within the literature (Williams, 1995, 1998; McGowan & Davis, 2019; McGowan & Tall, 1999) and labeled each concept within the ten maps using the features of mathematical discourse that make it distinct from other discourses as described by Sfard; as a reminder, these features are word use, symbolic mediators, narratives (inclusive of endorsed narratives), and routines. Routines are further categorized into explorations, deeds, and ritual. Explorations are further categorized into construction, substantiation, and recall. I restricted the initial analysis to the level of mathematical discourse features and the sub-level of routine types. Among the seven participants, two were instructors of calculus and the remaining five were students. Two students were enrolled in a collegiate Intermediate Algebra course and created two to three concept maps over the course of a semester. The remaining three students were a mixture of traditional and reform Calculus students and their concept maps were created at the end of the course. Propositions within the maps were labeled as narratives while the other features of the maps were labeled relatively easily in an ad hoc fashion. However, my decisions to label portions of a map as deed or ritual as opposed to exploration were not supported in the absence of interview data. I decided to further analyze the maps to generate ideas about ways to operationalize these labels in a manner similar to Morgan & Sfard (2016).

**REDACTED IMAGE**

Figure 8: Reform Calculus instructor's expert function concept map. I re-imagine Williams' (1998) example of a concept map as one composed of mathematical objects and the relationships between them. Reprinted from Williams, C. G. (1998). Using concept maps to assess conceptual knowledge of function. *Journal for Research in Mathematics Education. 29*(4), 414-421.

Due to my use of concept maps and realization trees outside of their intended structures, I sought support for applying the prior literature on the validity of concept maps to the discursive map. To accomplish this, I needed to establish that discursive maps would exhibit a strong similarity to the representations a participant would provide when asked to construct a concept map. I initiated this effort by studying the expert concept map for function in Figure 8 due to the feature of hierarchal organization (Williams, 1998). I began by making the basic assumption that unrestricted concept maps and discursive maps both aim to approximate the full potentiality of a phase, word, symbol, etc., subject to a specific academic field without invoking additional phrases that prompt one to act (i.e., a task-situation). I then considered the question of how a participant may attempt to represent the full potentiality of a phrase or word within their discourse without assistance. Williams (1998) and Novak & Cañas (2008) both specify asking concept map authors to begin by brainstorming a list in response to a signifier. They are then asked to organize this list into a concept map using linking words as necessary. Without a task-situation to organize this list, I reasoned that the map author will likely need to cycle through regularities in their list before they are able to organize their brainstormed realizations into a concept map. Therefore, I sought to identify regularities in Figure 8 using my function discourse.

    Without consulting commognitive constructs, I developed categories to encompass all features of the map using my perspective as an insider to expert mathematical discourse. Due to the inclination of categorization on the part of expert map authors, my initial categories largely followed the categories established by the map author. These were operations, common types, properties and definitions. I then categorized features of the map that were not explicitly named by the map author. The expert author of Figure 8 represented a cycle of terms around function that, based on the bi-directionality of links, appeared to be considered as equal or equivalent to the function signifier. The remaining map features were Fourier series, Taylor series, implicit form of a function, explicit form of a function, table is displayed as function, and function individually as graph is a common function. I reasoned that if the map contained a section for common types of functions, then any other functions that were not present in this category may have been considered uncommon, or, a general category of function that cannot be further categorized without being represented (i.e. implicit and explicit in the map encompass all yet-to-be defined functions). The expert map specifically contains the propositions: (1) Fourier series are a form of functions that are built from limits (also known as differentiate) which is an operation on function, and these series are linked to trigonometric functions which are common types of functions; and (2) Taylor series are a form of functions which are built from powers (which are built from polynomials) where polynomials are a common type of function. I identified common types, Fourier series, Taylor series, implicit, and explicit as instances/examples of function, however, the specification that the two series were *built from* another item in the map warranted further examination. My initial categories are below:

| Map regularities | How to identify |
|---|---|
| equivalent to (other names) | interchangeability with initial phrase or word indicated by equalizing language or bi-directionality at the same level of hierarchy |
| operations on (what can be done to it) | ways to manipulate the initial phrase or word, indicated by an action or process that happens "to" or "on" the signifier |
| properties/definitions (how it is talked about) | descriptors/ways to describe an instance of the signifier |
| Examples (how it "looks") | represents an instance of initial phrase or word, indicated by "is a", "are a", "type of", links and uni-directionality away from signifier toward example, subordination |

Figure 9: Preliminary categories to describe common features in concept maps.

    With the single expert map categorized, I then applied my initial categories to a reform student's function map (Williams, 1998, p. 417) and a traditional student's function map

REDACTED IMAGE

Figure 10: Traditional Calculus student's function concept map. Reprinted from Williams, C. G. (1998). Using concept maps to assess conceptual knowledge of function. *Journal for Research in Mathematics Education. 29*(4), 414-421.

(Williams, 1998, p. 418). Williams identified the traditional student's map as indicating the least conceptual understanding among all participants; this is inclusive of the expert group and both of the reform and traditional student groups. Unlike the expert's map, the students did not include series as functions. The traditional student did have "operations" as a linking word, though they only included add, subtract, multiply, and divide among function operations. I returned to my consideration of Fourier series and Taylor series in regard to their separation from common types of functions in the expert map. I sought to determine how to prod a student with Taylor series and Fourier series in their function discourse into representing this on their map without directly interfering with their map creation. This would seemingly require me to include instructions within the discursive map training that could encourage self-prodding on the part of the student every time that they create a map. I wanted to keep the discursive map instructions consistent between each map so that the only difference in map creation would lie in the signifier given.

  When examining the expert's map, I initially thought that it would be fitting to describe Fourier series and Taylor series as compositions of other functions. However, neither student appeared to view derivatives as a limit process, and therefore, asking students to consider compositions of their brainstormed list when grouping would be unlikely to yield the inclusion of, for example, *Taylor series* even if Taylor series were part of the student's discourse. Of note, the students did include *derivative* and *differentiation* in their respective maps. Using the traditional student's map, I determined that the operations of *addition, subtract, multiply*, and *divide* along with *differentiation* could possibly yield a Taylor series if this was in the student's discourse by considering my list of regularities dependent on the manner that differentiation is organized within the discourse. However, *differentiation* was represented in the student's map in a restricted manner. Namely, *differentiation* was only related to concavity, critical points, extrema, and partial differentiation; it would be difficult to coordinate the features of this student's map to create a Taylor series. However, this student exhibited a notion of first derivative and second derivative. This led me to pull away from the category *operations on*, and I instead replaced this with *interactions between instances.* Certainly, this would be inclusive of most operations and composition. However, this category still did not provide a path between derivative, first derivative, and second derivative as depicted in this student's map.

Using only the traditional student's map as evidence, I surmised that such a student could likely see derivatives as an action of a function on itself, not with another instance of itself, nor between itself and another function. This led me to further parse *interactions* into interactions between instances, actions within a single instance that create a new instance, and actions within an instance that preserve the current instance. Using the *instance* language, I arrived at my second framework of categories that are commonly present in concept maps:

| Map regularities |
|---|
| instances of word that are equivalent to word |
| the result of an interaction between instances of word |
| the result of an interaction between word and different word that was independently produced by map author |
| the result of an interaction between instances that creates a new instance of word |
| actions within an instance that create a new instance |
| actions within an instance that preserve the current instance |
| how an instance is talked about |
| instances that fall into the category of word |

Figure 11: Refined categories to describe common features in concept maps.

At this point, the regularities had become increasingly descriptive, and I removed the separate identification strategy from my framework. I was able to categorize most features of the maps using one of these regularities and all features of the expert instructor's map using one of these regularities. All three maps in Williams (1998) made mention of graphs in some fashion, though none of these maps included a graph. In anticipation of participants in my study including drawings in addition to symbols that they cannot define past knowing how to act in response to the symbol, I added a category for iconic representations of instances to my regularities. At this point, I worked through the remaining six concept maps and ensured that each element of the map could be labeled with one of the map regularities.

With the successful categorization of all ten concept maps, I returned to the commognitive framework and used the map regularities (relative to the parent signifier) to operationalize each commognitive construct. I came up with the following codes:

| Commognitive Construct & Source | | Map regularities |
|---|---|---|
| Author as Source | External Source | |
| (Interchangeable) Word use | (Interchangeable) Routine ritual | instances of word, drawing, symbol that are equivalent to word |
| Routine (exploration or deed) | Routine (deed or ritual) | the result of an interaction between instances of word |
| Routine (exploration or deed) | Routine (deed or ritual) | the result of an interaction between word and different word that was independently produced by map author |
| Routine (exploration or deed) | Routine (deed or ritual) | the result of an interaction between instances that creates a new instance of word |
| Routine (exploration or deed) | Routine (deed or ritual) | actions within an instance that create a new instance |
| Routine (exploration or deed) | Routine (deed or ritual) | actions within an instance that preserve the current instance |
| Narrative or endorsed narrative (cite endorser) | Normative endorsed narrative (if attributed to expert community) or routine ritual | how an instance is talked about |
| Routine (exploration or deed) | Routine (deed or ritual) | instances that fall into the category of word |
| Visual mediator | Visual mediator | isolated mathematical symbol |
| Visual mediator | Visual mediator | isolated drawing |

Figure 12: Refined categories to describe common features in concept maps and discursive maps.

Note, I developed these codes to label each mathematical object in a discursive map in relation to the initial signifier or a parent signifier. The map's resulting propositions form the map author's narratives on the discourse. The initial signifier is automatically labeled as *signifier* or *word use*. All drawings and symbols in isolation should be labeled as visual mediators. The remaining objects must be labeled as word use, routine, or narrative. The label assigned to the object corresponds to the justification the participant gives for the proposition.

     In the absence of supporting interview data, any target object that is considered a theorem, definition, axiom, etc., within normed mathematical discourse should be labeled as endorsed narrative. Next, any target represented in sentence or long form is labeled narrative. Finally, the remaining objects are labeled as routine without qualification. No determination is made regarding the routine sub-categorization. In the presence of supporting data, all target objects within one degree of separation from the initial signifier are categorized according to the regularity. The selection of a regularity should be supported by a participant's utterances. Then, dependent on the participant's description of the source for the proposition governing the target, the appropriate label under author as source or external authority as source should be selected. If a participant cites an external source as the reason for an element or narrative (inclusive of linked propositions) in their map, the source of the element is deemed external. A target mathematical object within a proposition that is the result of the participant seeking to sustain a bond with

others, or, if they are unable to further parse or provide an explanation for the object or proposition, the element should automatically be labeled as ritual. Realizations that can be further justified by the student should either be labeled exploration or deed at this level. Purely discursive routines should be labeled as explorations, while practical routines should be labeled as deeds (Sfard, 2008). Note, ritual is not an inherently negative characterization. Much of the mathematical discourse development in students begins as ritual (Sfard, 2008). If a participant realizes exponential as $e^x$ for no other reason than being taught that this is how an exponential function is represented (despite being false in general), that object should be labeled as ritual. As a final remark, my analysis of the expert instructors' discursive maps primarily featured word use, routines, and endorsed narratives as labels. Visual mediators were not present.

    The goal in analyzing student maps in this manner is to categorize targets of propositions that are one degree of separation away from the initial signifier. Of course, discourse is recursive, so in the same manner that close objects are categorized in my framework of codes relative to the initial signifier, the subsequent realization can again become a signifier. All close connections to this new signifier can be categorized in the same way. It is common for student function concept maps to unpack successive objects and branch further away from the topic of function (Williams, 1998). Experts tended to focus their maps on function specifically. The reader can compare Figures 8 and 10 to see an example of this. As a final verification that I had captured the typical features present in concept maps with my codes, I then applied my framework of codes to the algebra map in Figure 4. I was able to categorize all features of this map using the framework.

    According to Williams (1995; 1998) students consistently created concept maps that did not align with their initial training; this occurred despite confirmation that an initial concept map for *fraction* was produced correctly. As such, I expect that students will include connections to non-mathematical objects within the polygons of their discursive maps that fail to form a proposition that the student deems valid. For every non-mathematical object in the discursive map, if the student indicates that they intended to use the non-mathematical word in a polygon as a link between two mathematical objects, the proposition should be read as mathematical object – non-mathematical object – mathematical object provided that there are no linking terms on the lines forming the proposition. The target objects will still be considered as having one degree of separation from the initial signifier in this instance. Any propositions that fall outside of these specifications must automatically be labeled as a routine. If the participant cannot parse such a proposition so that it falls into the structure as described in training, that proposition must automatically be labeled as a ritual, and the participant should be queried about the source of that proposition within their discourse.

    With all commognitive constructs operationalized, I conclude by compiling a list of succinct, non-specialized instructions for inclusion in the discursive map training based on the final code list. These instructions are meant to encourage the participant to determine relationships between their list of realizations in response to the signifier. During the presentation on discursive maps, students will be told that discursive maps typically feature: (1) the most generalized categories of examples of the signifier to exhaustion; (2) terminology or descriptions that describe features, regularities, what is allowable, and what is not allowable regarding the signifier in general terms; and (3) what can be created by combining ideas from their brainstorming session. These combinations should include how the signifier interacts with other objects, how the signifier acts in relation to itself to create a new instance of the signifier, and how the signifier acts in relation to itself that leaves the signifier in an equivalent form. These instructions are given during the training and can be repeated to participants at any time in the

study if they are struggling to create either a discursive map or tree. With these specific instructions, I expect that discursive maps will yield the same data in the same representative structure as concept maps.

### Discursive maps in this study.

Students enrolled within an introductory proof course are offered tutoring in the form of problem sessions in exchange for their participation in Phase 1 of the study. As a result of conversations with mathematics education students that were previously enrolled in this type of course, I intend to target mathematics education students enrolled in Introduction to Higher Mathematics at the University of Georgia (UGA) for this phase. I will hold loosely-structured problem sessions within Aderhold Hall one to two times a week beginning around late-February and ending when the Spring semester concludes. Students may choose to attend or not attend the problem sessions as they deem necessary. Students that attend a problem session are required to come prepared with a task-situation from their course that they are struggling with, or in the absence of bringing their own task-situation, they will receive a task-situation based on a wider topic within the course that they have identified as a source of difficulty. Past and current syllabi for the version of this course that is taught at UGA have been obtained so that problem sessions can be planned according to the course schedule and focus.

Student participants in Phase 1 of the study will be asked to draw evoked trees to illustrate their initial thinking about their task-situation. This phase of the study allows for direct intervention and guidance for student participants while they construct proofs. Upon identification and interpretation of as many words or symbolic artifacts as student participants can realize within the task-situation, I will scaffold interventional instruction by having participants explore their realizations in response to my directed questioning of their progress. Students may be asked to create a discursive map in addition to the evoked tree if the participant is having difficulty beginning their tree or in the event that I need further information as to how this student has organized a mathematical object within their discourse.

During the timeframe where injective and surjective functions are covered in the students' transition to proof class, I will provide tasks similar to the tasks to be used in Phase 2 of the study for the Phase 1 participants to complete as practice. Any drafts of trees, discursive maps, or proof constructions created by Phase 1 participants will be reproduced and original copies of artifacts will be returned to the participant. After the conclusion of this phase of the study, maps, constructions, and trees produced in this phase are de-identified, assigned pseudonyms, and digitally re-created. Artifacts collected in this phase will be analyzed to identify common or typical student solution paths, word use, and routine features in their evoked discourse. This analysis will allow me to: (1) set expectations for the types of solutions and trees Phase 2 participants may produce; (2) add task-situations to the task bank to be used in Phase 2 that were sources of rich discourse reproduction by Phase 1 participants; and (3) identify areas of task interpretation and response difficulty within Phase 1 participants that may also arise in Phase 2 participants. Upon completion of the initial analysis, copies of artifacts that are relevant to function, injections, or surjections will be used to further de-identify student artifacts in the second phase of the study.

Phase 2 of the study begins with an on-campus discursive map training for consenting students and the course instructor. These trainings have been modeled according to Williams' (1998) description of concept map training of participants to produce traditional student, reform student, and expert function concept maps. The instructor will be trained prior to the course, while student participants will be trained according to their availability separate from the instructor. One mass training for students will occur toward the beginning of the course with subsequent student trainings being conducted if fewer than ten students consent to and are eligible for the study. I believe that beginning the class-based study with ten students will make concluding the study with at least three participants a realistic goal. During training, students are shown several images of different types of maps that have been called "concept maps" within education research across disciplines. These images will include spider maps, hierarchal concept maps, and non-hierarchal maps along with descriptions of each of these maps. Participants will be informed that the requirement for their maps is that a term, image, or symbol must be presented within enclosed polygons. These polygons can be connected to other polygons via lines. Students should only write "linking words" (e.g. is, are, equivalent to, equal to…) on these lines. The resulting two polygons connected by a link should form a proposition that the student deems valid. Participants will also be provided with the list of instructions as developed based on my framework of codes. In line with Williams (1998), students will be asked to make a discursive map for fractions. Students will be given progressive feedback to assist them in building discursive maps that are consistent with the training specifications. Students who produce discursive maps for fractions that meet the training guidelines will be invited to the initial interview and asked to create a discursive map for function at home. Students are instructed to produce their map from memory; they should not place anything on the map if they do not have any familiarity with that mathematical object. Further use of discursive maps occurs during the interviews in this study.

**Task-Based Interviews**

Much can be revealed about individual student thinking through the task-based interview (Goldin, 1997, 2000; Maher & Sigley, 2014). In this research design, the researcher creates a task for the participant with their specific research questions in mind. Like the structured interview (Patton, 2015), the researcher anticipates the responses of the participants and pre-plans how the task-based interview will proceed for those responses. While the task-based interview is popular in proof research (e.g. Guler & Dikici, 2014) the task-based interview is unnatural in that it removes the student from their normal learning environment, therefore, any interactions during the course of the interview will likely be influenced by the presence of the interviewer (Sfard, 2008). However, the task-based interview is beneficial in that it allows the interviewer to focus on individual differences in approach among participants and to further explore those differences when they occur.

    **Task-based interviews in this study.**

The goal of the initial student interviews (Appendix A) during Phase 2 are for all student participants in the class-based study to discuss their initial function discursive maps. I will ask participants to walk through their creation of their map, discuss their connections, and clarify bi-directional and uni-directional links in their maps. Due to the expected differences between student maps, the semi-structured nature of this interview will allow me to prompt students to explicate their rationale for links between the objects in their discursive maps and comment on the validity of propositions formed by two mathematical objects connected by a single link or two mathematical objects connected by a non-mathematical object within a polygon. If students state that a proposition is not valid in their discourse, the student participant is asked to describe the valid proposition they intended to express by identifying its mathematical objects and successive links within their map. Prior research using concept maps to study student function concept indicate that student propositions are commonly expressed as links between multiple polygons despite concept map training that explicates that true propositions should be formed by two polygons and a single link (Williams, 1998). I will use my insider perspective as an expert within the mathematical community to identify propositions that are not normed within expert mathematical discourse within the initial analysis. While I will not indicate this to participants, I will ask participants about the validity of these propositions within their discourse. The student participants' initial discursive maps are the closest approximation I can obtain of the students' function discourse prior to any instruction in the introductory proof course. This interview will help establish a baseline for the objects the student participants associate with function, thus allowing me to explore changes to these associations after instruction. Using the criterion of the strength of student explication of their rationale for their map, a subset of ten students invited for this initial interview will be invited to continue in the study throughout the semester. The instructor's discursive map will serve as the comparative map for students' function discourse; students' final maps will serve as the comparative map to students' initial maps.

     A central goal of the initial instructor interview (Appendix B) is to establish the instructors' function discourse in addition to determining the instructors' expectations of and for students in the course. Prior to this interview, the instructor will construct their personal discursive map for function, discursive maps that represent the function discourse that they believe is representative of students entering their introductory proof course, and a discursive map that they expect for their students' function discourse at the end of the semester. These artifacts will be collected, at the latest, one week prior to the initial interview with the instructor for analysis. During the interview, the instructor will discuss their personal map as well as the maps they expect of and for student function discourse at the beginning and at the end of the semester. I will also ask the instructor to discuss their experiences with student approaches to proof construction and student mastery of injections and surjections in particular. This interview will end with a discussion of the instructors' expectations for the development of students' proof discourses throughout the semester. During this time, the instructor will create a proof discursive map that they expect students to start with in the course.

     Task A is the focus of the second round of student interviews (Appendix C) and occurs after students have had some instruction in the introductory proof course, have covered set notation, and have yet to begin instruction about injections and surjections in this course. During this interview, I will ask students if there are any changes that they wish to make to their function map. If students produce any new maps, I will ask them to discuss their changes. Students that do not incorporate proof discourse into their function map will be asked if they can create a discursive map for proof. The purpose of refraining from asking students about proof prior to this interview is to allow students to fully explicate their function discourse while reducing interviewer influence on how that discourse is organized. This can aid in identifying if a student subsumes proof discourse into their function discourse, their proof discourse subsumes their function discourse, they connect their proof and function discourse but do not identify any one as subsumed by the other, or if students represent these discourses separately. I anticipate that students will either produce a separate proof map or incorporate function into their proof map. When a student completes their proof map, linked to their function discourse or not, I will inform students that they are not required to place both objects in one map, place one map within the

other, nor to keep them separate. I will then ask the student participant if they are satisfied with their map(s) as a representation of their current function and proof discourses. When the student is satisfied with their map(s), they will be asked to prove or disprove that a function, likely represented as a directed graph or arrow diagram (e.g. Dubinsky & Wilson, 2013), satisfies a given definition. Injection and surjection terminology is not used in this task, but the definitions of these are included in the task unlabeled. If the student participant completes the task, I will ask them if they can construct a tree for how they approached this task, though they may still make changes to their discursive map if they do not feel that their discursive map accurately captures their function discourse.

Students who are unable to begin proof constructions during Tasks A, B, and C will be asked to identify any key terms that they believe will assist them in interpreting the task, and to begin constructing a realization/symbolic tree. The tree construction can stop at any point where the student no longer has realizations to add or states they are unable to progress on their proof construction. Further participation of a student in this category will be determined on a case by case basis dependent on the potential of the students' incomplete tree, incomplete proof construction, and their discursive maps at this point of the study to produce interesting artifact triples for later examination by the course instructor. A student will not be removed from the study solely for being unable to construct a proof or tree for a single task. However, it is imperative that student participants can make initial progress on building a symbolic/realization tree through the identification and realization of words or symbols within the task, even if these realizations do not appear to contribute to a viable construction. The students' artifact triples are collected at the end of the interview, de-identified, re-constructed digitally to anonymize the participant, and labeled with a pseudonym. An initial analysis is conducted on all student triples for Task A prior to the instructor interview for Task A. Furthermore, the task for the next interview is selected from a pre-organized proof task bank upon this analysis (see Appendices C-E for representative samples from the task bank). The student interview structure, artifact triple creation and collection, and subsequent analyses for Task B and Task C largely follow the design of the Task A interview. Differences between subsequent student task-based interviews are described below.

Task B is the focus of the third round of student interviews (Appendix D). These interviews will be conducted after instruction on proving or disproving that a function is injective and surjective. In the third interview, students are presented with their most recent discursive map(s) and asked to make any changes that they feel are necessary. The student will then be presented with a task that requires them to prove that a function is injective and that another function is surjective. Students are not initially given definitions of these terms. To select this task, I collected twenty introductory proof course final exams from multiple universities and analyzed the types of injection and surjection problems that appeared on these exams. This type of task with minimal difficulty was prevalent in the sample final exams. The standard solution for this task involves completing a direct proof of injectivity by evaluating the given function at two arbitrary points and then deriving an equality. If the student asks for definitions of injective and surjective functions, they are given two definitions for each term. If the student participant completes the task, I will ask them if they can construct a tree for how they approached this task, though they may still make changes to their discursive concept map if they do not feel that their map accurately captures their proof and function discourses.

Task C is the focus of the fourth round of student interviews. These interviews will be conducted approximately one week prior to the end of the semester. The goal of this interview is to collect final discursive map(s) from students after a full semester of introductory proof and complete instruction on injective and surjective functions. Student participants are presented with their prior discursive map(s) and asked to make any changes that they see fit. In this interview, students are given one of two possible novel tasks. The first possible task is of minimal difficulty and was selected from a prior commognitive study of students' mathematical discourses on examination tasks in proof courses (Thoma & Nardi, 2018, p. 169)[5]. This task asks students to prove that a function with domain and codomain in the integers is injective, surjective, both or neither. Thoma and Nardi found that students, having only seen examples of this style of problem from the reals to the reals, did not consider the domain in their solution. However, domain considerations are of critical importance to normative proof constructions in response to this task. Thoma & Nardi (2018; 2020) found that students relied on *precedent events* for proof construction in response to this task, and therefore, students were unable to complete this task to the satisfaction of their course instructor due to applying rituals for classifying functions as injections and surjections where both the domain and range were on the real numbers. This task will be used if student participants have exhibited substantial difficulty in constructing proofs in response to Task A and Task B.

      The second possible task is of significant complexity, has increased student autonomy from prior tasks, and was adapted from a previous study of student function development by Breidenbach, Dubinsky, Hawks, & Nichols (1992, p. 283). This task was selected due to its partial symbolic exploitability as well as requiring students to go beyond their symbolic approach to complete a normative proof construction. Morgan & Sfard define the grain size of a task as "the minimal number of decisions (choices) the problem solver must make while designing a series of elementary steps necessary to solve the problem" (Morgan & Sfard, 2016, p. 113). Elementary steps are those than cannot be parsed in a recursive manner into pieces where further decisions must be made on the part of the student. The increased grain size of this problem when compared to prior tasks in this study and student participants' likely unfamiliarity with a similar task should generate meaningful discourse exploration by student participants. Students are not initially given definitions in either task choice. If the student participant completes the task, I will ask them if they can construct a tree for how they approached this task, though they may still make changes to their discursive map if they do not feel that their map accurately captures their function and proof discourse. Due to the difficulty of this task, I expect that most students will need to create a tree to explore their discourse and engage in active decision making prior to writing a proof construction. This is in opposition to the previous tasks where I expect several students to automatically begin their proof construction based on *precedent events* within the proof course. As an alternative and dependent on participant ease with Task A and Task B, I may include the Thoma & Nardi option as a second problem within Task B, while assigning the Breidenbach et. al. task to all participants for Task C.

      The second, third, and fourth rounds of interviews have an instructor component and each of these interviews occurs between one and two weeks after the student interviews for one task situation conclude. Prior to the instructor interviews for Tasks A, B, and C, the instructor is asked to respond to each task situation with a proof construction that they believe will be typical of their students. Instructors can create up to three different proof constructions that they believe will be the constructions most similar to their students' approaches. For each construction, the instructor creates a tree that includes the key words the instructor expects their students to use to

---

[5] A similar version of this task can be found in Appendix C. Images of tasks that have been selected from prior studies are labeled and included in the specified interview protocols within the appendices.

interpret the task and construct a proof, along with successive realizations until the proof can be constructed. Each realization tree should represent the instructor's best approximations of how students in their transition to proof course are likely to approach the task situation. I will collect the instructor's proof construction and tree pairs for analysis at least two days prior to the instructor interview for the associated task. This delay will allow me to analyze the elements of mathematical discourse present, using my insider perspective, within the instructor's pair of artifacts for this task. These artifacts will be compared to the student participants' artifacts for further analysis and for the selection of student artifact triples to be shown to the instructor during this interview.

    The one-two week time delay between student and instructor task interviews for a single task will allow me to conduct an initial analysis of the student participants' artifact triples. As a reminder, the artifact triples for each task consist of the student's proof construction in response to the task, the associated tree, and the discursive map used by the student at the time of the interview. Only the proof construction portion of a student's artifact triple will be selected for further review by the instructor within these three follow-up interviews if, by my insider perspective as a member of the expert mathematical community, I determine that elements of the student's discourse are not consistently represented between all three artifacts created by the student. Additional flexible (i.e. not required) considerations in the selection of student triples are made based on the strength of: (1) the student participant's likelihood of participating in the remainder of the study so that a complete sequence of maps can be analyzed and discourse expansion can be examined over time in fulfillment of first research question; (2) the student participant's explication and representation of their discourses as evidenced by their prior task-based interviews so that I do not violate of the principles of commognition by attempting to infer student reasoning; (3) the alignment of the students' proof constructions with normed constructions or the instructor's construction while simultaneously demonstrating discursive maps that represent highly restrictive and ritualized discourses of function and proof; and (4) sparse or incomplete proof constructions while simultaneously demonstrating discursive maps that are rich, detailed, or normative. I expect that the number of constructions presented to the instructor will gradually decrease across all tasks due to a reduction in student participants over time.

    During these follow-up interviews with the instructor (Appendix F) for each task, the instructor is asked to discuss their proof construction(s) and their expected student trees. This interview will allow me to clarify the discourses that the instructor(s) expect students to draw on to construct a proof in response to the task. The protocol for this interview emphasizes the elements of mathematical discourse (i.e. word use, visual mediators, narratives, and routines) that the instructor expects their students to draw on to interpret the task situation. In the latter half of this interview, the instructor is provided with de-identified, digitally reconstructed proofs from pre-selected student artifact triples of their students, Phase 1 participants, and volunteers. I will ask graduate students and math educators at the University of Georgia to contribute artifact triples for this purpose. The instructor is asked to examine each proof construction and to verbally describe the features of these constructions that the instructor is using to construct a narrative about whether this student has a function discourse and proof discourse that are, as an example: (1) exceeding expectations; (2) meeting expectations; (3) developing; (4) not meeting expectations as defined by the course instructor. The narrative ratings used by the instructor to determine artifact acceptability are negotiated with the instructor prior to the beginning of the course. If the instructor does not want to choose their own system of narrative ratings, the

aforementioned example narratives will be used by the instructor. The instructor will be asked to describe how they think the student arrived at their construction. The instructor will then finalize their narrative about this students' discourse by labeling the proof construction with this narrative. The marked proof constructions are then returned to their associated student triple for further analysis when this interview concludes.

After the conclusion of student and instructor interviews for Task C, student artifact triples to be presented to the instructor in the final semi-structured clinal interview will be selected using the following flexible (i.e. not required) criteria: (1) the student participant's participation in all four interviews; (2) the student participant's ease of explication and representation of their discourses as evidenced by all four of their interviews; (3) the student participant's completion of all three proof constructions without regard to accuracy; (4) the likelihood, based on my insider perspective as a member of the expert mathematical community, that the instructor's narratives on each artifact in the student's artifact triple will differ; (5) minimal risk of instructor identification of the student by examining their artifacts. If greater than five student participants are retained in the study through the end of Task C, the number of triples to be presented to the instructor will be selected dependent on the length of time the instructor is willing to dedicate to the final interview. The final interview will be at least one hour. If less than or exactly five student participants are retained in the study by the end of the semester, all remaining student participants will have their artifact triples shown to the instructor in the final interview.

The final interview (Appendix G) will be conducted after the semester concludes and will only be held with the instructor. Only selecting from student artifact triples that have proof constructions labeled with the instructor's narrative rating, the instructor is shown the discursive map used by the student at the time of the associated proof construction that the instructor has assessed. These discursive maps are presented to the instructor in an order randomly chosen and unrelated to the order in which they were shown student proof constructions. The instructor will be asked to examine the discursive concept map of the student and decide if they believe that a student with this map would have been able produce an acceptable proof construction to the related task, and if this student's discourse is what the instructor expects of a student at the time in the course when the map was constructed. The instructor is asked to verbally discuss their reasoning using evidence from the student's discursive map and the task situation. Using the same narrative labels, the instructor marks the discursive map with one of the pre-negotiated narrative labels. This is done for all of the discursive maps associated with the marked proof constructions selected for this final interview.

The instructor is then presented with a marked proof construction and allowed to re-read the construction and their narrative rating to re-familiarize themselves with their rating. The instructor will then be presented with the associated student discursive map from this student's artifact triple and asked for commentary regarding their narratives on each related artifact. As a reminder, the instructor has rated this discursive map during this interview. Finally, the instructor is provided with the student tree from this artifact triple and asked for commentary. This procedure is repeated for all student artifact triples that have been selected for inclusion in the final interview. The goal of this interview is to document the instructor's responses to how well their assessment of a students' function and proof discourses based on their reading of a students' proof construction aligns with their assessment of these discourses using the student's discursive map. Furthermore, I would like to document the instructor's commentary when presented with the tree from this students' artifact triple.

**Data Collection**

  Phase 1 participants are asked to draw restricted trees in response to participant selected task-situations as well as proof constructions. These participants are additionally given injection and surjection proof task-situations and asked to create trees and proof-constructions in response. At this phase, participants may be asked to create discursive maps. Any drafts of trees, discursive maps, or proof constructions will be reproduced and original copies of artifacts will be returned to the participant. After the conclusion of this phase of the study, maps, constructions and trees are de-identified, assigned pseudonyms, and digitally re-created. Upon completion of the initial analysis, copies of artifacts that are relevant to function, injections, or surjections will be used to further de-identify student artifacts in the second phase of the study.

  Documentation of the classroom environment and examining artifacts are helpful in gaining insight into student discourse development. During Phases 2 and 3 of the study, I will take detailed field notes on the classroom environment, lecture material, student assignments, student activities, and course texts to document the classroom microculture. These notes, combined with the video recording and transcription of function lectures, interviews (whether semi-structured clinical or task-based), instructor-produced artifacts, and examining student-produced artifacts from student interviews, provide a rich set of data from which to triangulate and support theorized contributors to student discourse development. These data can also provide context for utterances made during one-on-one interviews with students and instructors.

  The classroom environment and lecture will be recorded on each day that the function unit is taught. To best capture the interactions and activity in the classroom, I will use two cameras in the classroom to provide different views of the environment and the conversations that occur within it. The instructor for the course will be interviewed before the course, at three points during the semester, and after the course. The student participants will be interviewed at the beginning of the course and at three points during the semester. I will act as the participant interviewer during all interviews. While I will not attempt to alter student discourses, I describe myself as a participant interviewer due to the impact that interviewers can potentially have on the participant interviewee. Furthermore, the principles of commognition dictate that my utterances during interviews must be incorporated into the analysis (Sfard, 2012). Any interaction that students have with me during interviews and while they complete tasks are discourse developers, so these interactions must be accounted for. Each interview will be video-recorded and transcribed.

  During student interviews, students are asked to identify key words or symbols mediating their interpretation of the task situation and to construct symbolic/realization trees using these key words. I will use an insider perspective to determine narratives in a student participant's proof construction that have not been described by the student within their task-specific symbolic or realization tree. I will use an outsider's perspective to probe students and have them further explicate the relationships between successive mathematical objects within their trees and how students coordinate these objects to construct a proof. I will use a line of questioning with students that encourages their explication of their applicability conditions for their chosen routines, their course of action, and the conditions of closure.

  Prior to the second, third, and fourth instructor interview, instructors are asked to construct the proofs they expect their students to create in response to a task-situation. They are then asked to identify key words they believe students will use to interpret the task situation. The instructor will construct a symbolic/realization tree that they believe represents the likely signifier-realizations students will use to construct their proof. On the day of the instructor interview, the instructor will discuss their proof construction and describe the approach they expect students to use using their expected student trees. I will use a line of questioning with the instructor that encourages their explication of their applicability conditions for their expected routines, their expected courses of action, and the expected conditions of closure that students may appeal to complete their construction. Each instructor interview that follows student task-based interviews concludes with the instructor placing ratings, or narratives, concerning students' discourse development relative to the goals of the course.

The instructor's proof constructions and expected student trees in response to the tasks given to student participants will be collected and photocopied three times during the semester. These constructions will be typed, de-identified, and placed in a computer folder bearing the pseudonym of the instructor participant. Proof constructions, discursive maps, and trees produced by each student participant will be de-identified, photocopied, digitally re-created, and placed in the de-identified folder for each participant.

## Data Analysis

### Analyses of student data.

Phase 1 artifacts will be analyzed to identify common or typical student solution paths, word use, and routine features in their evoked discourse. This analysis will allow me to: (1) set expectations for the types of solutions and trees Phase 2 participants may produce; (2) add task-situations to the task bank to be used in Phase 2 that were sources of rich discourse reproduction by Phase 1 participants; and (3) identify areas of task interpretation and response difficulty within Phase 1 participants that may also arise in Phase 2 participants. Upon completion of the initial analysis, copies of artifacts that are relevant to function, injections, or surjections will be used to further de-identify student artifacts in the second phase of the study. The remainder of this section refers primarily to the class-based portion of the study and any aspects of the Phase 1 study for which similar data is produced.

In line with the five principles of commognition, selected data will initially be presented without my narratives imposed on the data. Maps will be examined for the validity of propositions by the student participant. This will occur within the student's task-based interviews. Coherence of maps are therefore determined by participants. Initial and final student maps and instructor maps will be presented without interpretive commentary initially, though examinations of commonalities of mathematical objects within student maps will be described. The majority of the relevant analysis is completed in the final interview by the instructor of the course. The instructors' narratives on student artifact triples are initially presented without interpretation.

I will identify the mathematical objects and propositions present in the discursive concept maps of each student participant to compare student discourses early in the semester so that I can establish changes in an individual student's discourse between the beginning and the end of the semester. For the three time points where students are interviewed during the course, I will visually represent whether a participant made a change to their map and the nature of the change. I will determine the number of students that made changes and the typical nature of those changes at each of these time points. Key word identification within task situations and associated realizations will be compared across student participants. To maintain the integrity of student maps and their intended discourse representation, I will use videos and transcripts of student interviews with additional support from any student interview where the student participant further explicates how they constructed their initial map, their reasoning for making insertions, deletions, links, and/or structural reorganization to their discursive maps throughout the semester. I will further support the identification of the student's discourses by using any series of utterances between the participant and interviewer where students comment on the validity of propositions formed by two mathematical objects connected by a single link, or a proposition formed by a non-mathematical object within a polygon between two mathematical objects, as well as student-identified propositions formed by successive links between more than two mathematical objects.

      Using a constant comparison approach, I will develop codes for similar objects, propositions, propositions deemed valid by students, and linkage present in student maps to further describe similarities and differences in function discourse among and between students at the beginning of the course. I will coordinate classroom field notes, syllabi, instructional resources offered to students, video recordings and transcripts of function lectures, and the video and transcripts of student interviews to develop categories of student experiences between interviews that may have resulted in map changes. This approach will also be used to suggest common linguistic relationships used among students to engage in proof construction for each task. Suggested theories will be supported with video and transcripts from student task-based interviews, classroom field notes, student discursive maps, student symbolic/realization trees, student proof constructions, and video recordings of function lectures.

      The video and transcript of the associated task-based interview will be used to classify each portion of student trees into categories of routine. If students are unable to construct partial symbolic/realization trees that consist of successive object realization branches that represent their proof construction strategy, this branch will be classified as either a deed or ritual. Students who appeal to explanations that specify that they are engaging in mimicking behavior to explain branches in their symbolic/realization tree will also have those branches determined to be ritual. Support for or against the appearance of different routines in student artifacts will be gleaned from classroom field notes, instructional resources provided in the class, and transcripts of interviews where students describe the resource that results in an identified branch of their symbolic/realization tree or discursive map. Instances of routine will be quantified. By triangulating student descriptions in transcripts, proof constructions, classroom field notes, and student discursive maps, I will suggest the source, whether instructional, student selected, or otherwise of each instance of routine. The quantity, type, and source of each routine will be compared across students for each task situation to suggest commonalities in routine selection between students.

      **Analyses of instructor data.**

      I will identify the mathematical objects and propositions present in the instructor's discursive map, expected student trees, the initial expected student discursive map, and final expected student discursive map. To maintain the integrity of the instructor's map(s) and their intended discourse representation, I will use videos and transcripts of the instructor's initial interview with additional support from any instructor interview where the instructor participant further explicates how they constructed their map and expected student trees. I will support the identification of the instructor function discourse by using any series of utterances between the participant and interviewer where an instructor explicates their rationale for links between the

objects in their discursive maps and comments on the validity of propositions formed by two mathematical objects connected by a single link, as well as instructor-identified propositions formed by successive links between more than two mathematical objects. Using a constant comparison approach, I will propose categories of objects, propositions, and linkages present within the instructors' discursive map of function.

The analysis of the initial and final expected student maps is to establish changes the instructor expects in students' function discourse between the beginning and the end of the semester. I will represent the instructors' expected student map changes and the suspected nature of these changes in a chart. I will support my determination of the nature of these changes and the course instruction that the instructor expects to result in these changes using video and transcripts of interviews where the instructors makes these goals explicit and explicates their reasoning for making insertions, deletions, and other linkage and/or structural reorganization to their initial expected student map to arrive at the final expected student map. The instructors' stated goals will be compared to the intended course schedule as represented by the course syllabi and the realized course schedule to suggest timelines of function discourse alteration expected by the instructor.

I will use the expected student discursive maps, expected student proof constructions, and expected student trees to identify mathematical objects and propositions that the instructor expects their students to appeal to in response to the task-situation. I will compare mathematical objects and propositions that are present in the instructors' expected student trees across tasks to suggest the discourses that the instructor finds necessary to construct a proof in response to the task situation, and the combined discourse necessary to construct a proof in response to all three task situations.

The video and transcript of the instructors' task-based interviews for Task A, Task B, and Task C will be used to classify each portion of the instructor's expected student trees into categories of mathematical discourse. These are word use, narrative, symbolic mediators, and routine. Though unexpected, if the instructor is unable to construct trees that consist of successive object realization branches that represent their expected student proof construction strategy, this branch will be classified as either a routine deed or routine ritual. Instructors who appeal to explanations that indicate that they expect students to engage in mimicking behavior to explain the branches in the expected student tree will also have those branches determined to be ritual. Support for or against using the four elements of mathematical discourse (word use, narrative, symbolic mediators, and routine) and their sub-elements to classify branches of the instructor's expected student maps, constructions, and trees will be gleaned from classroom field notes, instructional resources provided in the class, and transcripts of interviews where the instructor describes the instruction that they expect to result in an identified branch of their expected student tree. Instances of routine will be quantified. By triangulating instructor descriptions in interview transcripts, expected student proof constructions, classroom field notes, and expected student discursive maps, I will suggest the portion of instruction that the instructor likely expects to be the source of their students' use of specific routines. The quantity and source of each routine will be compared to student's quantity and sources of routines for each task situation to suggest commonalities and differences in routine selection and the source of routines between students and the instructor in the inter-group analyses of the data.

**Inter-group analyses of data.**

Mathematical objects and propositions common to the instructor's discursive map and students' initial discursive maps will be compared using the codes and categories created from the prior constant comparison approach to the analysis of the instructor and student data. In the same manner, the instructor's expected discursive maps for students and expected student trees will be compared to student discursive maps and student created trees for each task. Furthermore, mathematical objects and propositions that the instructor identifies as likely present or missing within a students' discourse based on their reading of the student's proof construction is compared to the actual student discourse map that was in use at the time of the proof construction.

I will produce a table representing the instructors' narratives about student discourses based on student proof constructions and student discursive maps. As a reminder, the instructor selects the narratives that they would like to use to rate student artifacts toward the beginning of the study. I will quantify the instructor's consistent vs. inconsistent narratives between artifacts based on differences in ratings provided by the instructor. I will produce a second table representing the instructor's narratives between the three artifacts and quantify consistent vs. inconsistent narratives between artifacts.

**Summary of Data Collection**

In summary, the proposed study will proceed in the following manner:

- Initial phase of study occurs from mid-semester until the semester concludes. Participant artifacts are analyzed to set expectations for the class-based study and to determine common solution paths and routines used by students.
- Phase two of study begins shortly before/at the beginning of the summer semester. Discursive map trainings with consenting students and course instructor. Function discursive map created by participants and shared with me. The course instructor selects narratives for artifact rating according to their preference.
- Videotaping of individual clinical interviews with all consenting students and instructor regarding their initial discursive maps. Clarify any directionality. Collect any new discursive concept maps that are created. The instructor also provides expected initial and final student discursive maps.
- Initial course instruction and creation of field notes. Student introduction to set theory by course instructor.
- Videotaping and transcription of second round of interviews conducted with each participant. Injection and surjection terminology are not used in the initial task, but the accompanying definitions are included in the task unnamed. Collect any new discursive maps that are created. Collect proof constructions and trees from instructor and all participants. Collect marked student proof constructions from instructor.
- Observation of course lecture and creation of field notes continues.
- Videotaping and transcription of third round of interviews on proving a function is injective or surjective. The function in this task should map from the real numbers to the real numbers and should be able to be completed by symbolic manipulation. Students are not initially given definitions for

- injective and surjective. If student asks for these definitions, they are given two definitions for each term. Collect any new discursive maps that are created. Collect proof constructions and trees from instructor and all participants. Collect marked student proof constructions from instructor.
- Videotaping and transcription of instruction on functions and proving/disproving functions are injective and surjective.
- One week before the end of the semester, videotaping and transcription of fourth round of interviews on novel task. If participants demonstrate substantial difficulty completing Task B, Task C should focus on a prove/disprove task that requires domain consideration (e.g. Thoma & Nardi, 2020). Otherwise, Task C should be of increased complexity (e.g. Breidenbach et. al., 1992). Collect any new discursive maps that are created. Collect proof constructions and trees from instructor and all participants. Collect marked student proof constructions from instructor.
- Semester concludes. The instructor submits students' course grades to their university and Phase 3 of the study begins. When grades are released to students, the instructor can be interviewed.
- Videotaping and transcription of final interview with instructor. Collect instructors' marked commentary on de-identified student artifact triads.
- (Optional) Phase 4 of the study begins shortly before the Fall semester. This phase repeats the class-based phase of the study in the event that a typical course instructor cannot be located, if protocols need to be refined, if it is deemed necessary, or in the event of a COVID crisis or other unforeseen circumstance. If the university of study needs to go online, a new IRB will be submitted to repeat the study online. However, Phase 4 will no longer be tied to a specific instructor under these circumstances.

**Summary of Study**

I reiterate that studies of injective and surjective functions in the context of proof at the postsecondary level are scarce. The dissertation of Thoma (2018) makes the most advances in the direction of the aims of my research by determining the characteristics of examination tasks among first year mathematics students, lecturer's expectations for students' engagement with university mathematics discourse and their use of these expectations in creating examination tasks, and the differences between students' school mathematics discourse and university mathematics discourse as well as the commognitive conflicts that could be observed as resulting from those differences.

While Thoma does not specify non-proof content to be the focus of study, all of Thoma's tasks require students to engage in proof construction. My study focuses on function, injections, and surjections in the context of proof. The consistent focus on injective and surjective functions allows for the study of students' temporal development of their function discourse in addition to the gradual incorporation of proof into students' mathematical discourse. This includes an accounting of the instructor and students' unrestricted function discourse toward the beginning of the semester, and the students' unrestricted proof discourse one month after the course begins. The commognitive view acknowledges that specific discourses are evoked in task situations, so it is essential that I collect students' expressed unrestricted discourses to obtain the closest approximation of the students' initial discourses as possible.

    The instructors in Thoma's study provided sample constructions and discussed their expectations of students in a semi-structured interview. Thoma appeals to the instructor's expectations for student proof constructions and the transcripts of interviews with instructors to identify the instructor's expected discourse for students. In my study, the instructor's expected proof constructions are supported by the instructor's hierarchical realization trees for the expected student approach to the task. These realization trees are provided prior to the interview, and therefore they are not influenced by questioning within my instructor interview protocols. I further support the determination of the instructor's expectations for students' discourse using transcripts of semi-structured instructor interviews. Thoma's proof construction and interview transcript strategy for the instructor is the same as their strategy for student participants. I additionally support students' representations of their discourses using the realization trees created by student participants. Namely, the product of the act of proof construction by students is referenced against the students' expressed strategy via their realization trees as well as their unrestricted discursive maps. This allows me to compare discursive features among all three student artifacts, including the discourse at rest and the evoked discourse. Furthermore, this allows for a comparison of the instructor's narratives regarding the student's development of their function and proof discourses between student artifacts. This feature of my study can potentially demonstrate the assessment advantages/disadvantages of discursive maps, realization/symbolic trees, and completed proof constructions in comparison to each other.

    The critical difference between our studies is that I have designed my study to produce data on students' discourses at three levels: (1) unrestricted discourse; (2) evoked discourse; and (3) implied discourse of the construction. This will allow me to comment on students' paths to proof construction and the extent to which they rely on their unrestricted discourses, if at all, to respond to task situations. When compared to the instructor's: (1) expectations for students' unrestricted discourse; (2) expectations for students' evoked discourse; and (3) expectations for students' unrestricted discourse based solely on a student's proof construction, I can compose narratives about the state of students' engagement with proof construction and developing non-proof mathematical discourses while transitioning to higher level mathematics.

| Research Question | Sub-Question | Data Collection | Data Analysis |
|---|---|---|---|
| In what ways do undergraduate mathematics students interpret and make sense of injection, surjection, and proof instruction in transition to proof courses? | How do students think about injections and surjections in the context of function before participating in transition to proof courses as evidenced by their initial discursive maps? | Initial discursive maps of students | I will identify the mathematical objects and propositions present in the initial discursive concept maps of each student participant to compare student discourses early in the semester. To maintain the integrity of student maps and their intended discourse representation, I will use videos and transcripts of the student's initial interview with additional support from any student interview where the student participant further explicates how they constructed their initial map.<br><br>I will support the identification of the student author's initial discourse by using any series of utterances between the participant and interviewer where students explicate their rationale for links between the objects in their discursive concept maps and comment on the validity of propositions formed by two mathematical objects connected by a single link, as well as student-identified propositions formed by successive links between more than two mathematical objects. Using a constant comparison approach, I will develop codes for similar objects, propositions, propositions deemed valid by students, and linkage present in student maps to further describe similarities and differences in function discourse among and between students at the beginning of the course. |
| | | Final discursive maps of students | |
| | | Video recording of initial student interviews and interview where final map is created | |
| | | Transcripts of initial student interviews and interview where final map is created | |
| | | Classroom field notes | |
| | What changes do students make to their discursive maps pertaining to proof, injections, and surjections as a result of lecture, assignments, and reading of the course textbook as well | Classroom field notes, notes taken on instructional resources (e.g. syllabi, textbook) and alternative resources used by students (e.g. textbook, digital content) | I will identify the mathematical objects and propositions present in the discursive concept maps of each student participant to compare student discourses early in the semester so that I can establish changes in an individual |
| | as the use of other student-selected resources over time? | All discursive concept maps that students create during study | student's discourse between the beginning and the end of the semester.

For the three time points where students are interviewed during the course, I will visually represent whether a participant made a change to their map and the nature of the change. I will determine the number of students that made changes and the typical nature of those changes at each of these time points.

To maintain the integrity of student maps and their intended discourse representation, I will use videos and transcripts of student interviews with additional support from any student interview where the student participant further explicates how they constructed their initial map, their reasoning for making insertions, deletions, links, and/or structural reorganization to their discursive maps throughout the semester. I will coordinate classroom field notes, syllabi, instructional resources offered to students, and the video and transcripts of student interviews to develop categories of student experiences between interviews that may have resulted in map changes. |
| | | Video recording of initial student interviews and interviews where discursive concept maps are altered | |
| | | Transcripts of initial student interviews and interviews where discursive concept maps are altered | |
| | How do students draw upon their injection, surjection, function, and proof discourses to interpret proof task situations involving injections and surjections? | Classroom field notes, notes taken on instructional resources, and alternative resources used by students | During student interviews, students are asked to identify key words or symbols mediating their interpretation of the task situation and to construct symbolic/realization trees using these key words. I will use an insider perspective to determine narratives in a student participant's proof construction that have not been described by the student within their task-specific symbolic or realization tree. I will use an outsider's perspective to probe students and have them further explicate the relationships between successive |
| | | Student proof constructions | |
| | | Discursive concept maps created by students | |
| | | Trees created by students | mathematical objects within their trees and how students coordinate these objects to construct a proof. I will use a line of questioning with students that encourages their explication of their applicability conditions for their chosen routines, their course of action, and the conditions of closure.<br><br>A constant comparison approach will be used to suggest common linguistic relationships used among students to engage in proof construction for each task. Suggested theories will be supported with video and transcripts from student task-based interviews, classroom field notes, student discursive maps, student symbolic/realization trees, student proof constructions, and video recordings of function lectures. |
| | | Video recordings of all student interviews | |
| | | Transcripts of all student interviews | |
| | | Phase 1 participant's trees, proof constructions, and discursive maps | |
| | What routines do students draw upon to construct proofs in response to task situations? | Classroom field notes, notes taken on instructional resources, and alternative resources used by students | Phase 1 artifacts will be analyzed to identify common or typical student solution paths, word use, and routine features in their evoked discourse.<br><br>During student interviews in Phase 2, students are asked to identify key words or symbols mediating their interpretation of the task situation and to construct symbolic/realization trees using these key words. I will use an insider perspective to determine narratives in a student participant's proof construction that have not been described by the student within their task-specific tree. I will use an outsider's perspective to probe students and have them further explicate the relationships between successive mathematical |
| | | Present discursive concept maps of students | |
| | | Student tree for present task | |
| | | Student proof constructions | |
| | | Video recordings of all student interviews | |
| | | Transcripts of all student interviews | objects within their trees and how students coordinate these objects to construct a proof. I will use a line of questioning with students that encourages their explication of their applicability conditions for their chosen routines, their course of action, and the conditions of closure. |
| | | Phase 1 participant's trees, proof constructions, and discursive maps | The video and transcript of the associated task-based interview will be used to classify each portion of student trees into categories of routine. If students are unable to construct partial symbolic/realization trees that consist of successive object realization branches that represent their proof construction strategy, this branch will be classified as either a deed or ritual. Students who appeal to explanations that specify that they are engaging in mimicking behavior to explain branches in their symbolic/realization tree will also have those branches determined to be ritual. Support for or against the appearance of different routines in student artifacts will be gleaned from classroom field notes, instructional resources provided in the class, and transcripts of interviews where students describe the resource that results in an identified branch of their symbolic/realization tree or discursive map. Instances of routine will be quantified. By triangulating student descriptions in transcripts, proof constructions, classroom field notes, and student discursive maps, I will suggest the source, whether instructional, student selected, or otherwise of each instance of routine. The quantity, type, and source of each routine will be compared across students for each task situation to suggest commonalities in routine selection between students. |
| In what ways do instructors | How do instructors think about injections and | Discursive map of instructor(s) | I will identify the mathematical objects and propositions present in the instructor's discursive map. To maintain the |
| interpret and make sense of proof, injective functions, and surjective functions? | surjections within the context of function as evidenced by their discursive maps? | All expected student trees created by instructor | integrity of the instructor's map(s) and their intended discourse representation, I will use videos and transcripts of the instructor's initial interview with additional support from any instructor interview where the instructor participant further explicates how they constructed their map and expected student trees.<br><br>I will support the identification of the instructor function discourse by using any series of utterances between the participant and interviewer where an instructor explicates their rationale for links between the objects in their discursive maps and comments on the validity of propositions formed by two mathematical objects connected by a single link, as well as instructor-identified propositions formed by successive links between more than two mathematical objects. Using a constant comparison approach, I will propose categories of objects, propositions, and linkages present within the instructors' discursive map of function. |
| | | Video recordings of all instructor interviews | |
| | | Transcripts of all instructor interviews | |
| | How do instructors expect students' injection and surjection discourses to be organized within the context of their function discursive maps prior to the transition to proof course? Including proof discourse, how do | Instructor's expected student initial discursive maps | I will identify the mathematical objects and propositions present in the initial and final expected student maps to establish changes the instructor expects in students' function discourse between the beginning and the end of the semester. I will support the identification of changes between the instructor's expected student maps using video and transcripts of instructor interviews where the course instructor explicates their reasoning for making insertions, |
| | | Instructor's expected student final discursive maps | |
| | | Video recordings of all instructor interviews | |
| | instructors expect students' discourses to be organized by the end of the transition to proof course? | Transcripts of all instructor interviews | deletions, and other linkage and/or structural reorganization to their initial expected student map to arrive at the final expected student map. |
| | | Classroom field notes and notes taken on instructional resources | I will represent the instructors' expected student map changes and the suspected nature of these changes in a chart. I will support my determination of the nature of these changes and the course instruction that the instructor expects to result in these changes using video and transcripts of interviews where the instructors makes these goals explicit and explicates their reasoning for making insertions, deletions, and other linkage and/or structural reorganization to their initial expected student map to arrive at the final expected student map. The instructors' stated goals will be compared to the intended course schedule as represented by the course syllabi and the realized course schedule to suggest timelines of function discourse alteration expected by the instructor. |
| | What discourses do instructors expect students to draw upon to interpret task situations involving injections and surjections in order to engage in proof construction as evidenced by their expected student trees? | Instructor's expected student discursive maps | Prior to their second, third, and fourth interview, instructors are asked to construct the proofs they expect their students to create in response to a task situation. They then identify key words they believe students will use to interpret the task situation. The instructor will construct a tree that they believe represents the likely realizations students will use to construct their proof. On the day of the instructor interview, the instructor will discuss their proof construction and describe the approach they expect students to use using their expected student trees.

The video and transcript of the instructors' task-based interviews for Task A, Task B, and Task C will be used to |
| | | Instructor's expected proof constructions for Tasks A, B, and C | |
| | | Instructor's expected student symbolic trees for Tasks A, B, and C | |
| | | Video recordings of all instructor interviews | |

| | | Transcripts of all instructor interviews | classify each portion of the instructor's expected student trees into categories of mathematical discourse. These are word use, narrative, symbolic mediators, and routine. Though unexpected, if the instructor is unable to construct trees that consist of successive object realization branches that represent their expected student proof construction strategy, this branch will be classified as either a routine deed or routine ritual. Instructors who appeal to explanations that indicate that they expect students to engage in mimicking behavior to explain the branches in the expected student tree will also have those branches determined to be ritual. Support for or against using the four elements of mathematical discourse (word use, narrative, symbolic mediators, and routine) and their sub-elements to classify branches of the instructor's expected student maps, constructions, and trees will be gleaned from classroom field notes, instructional resources provided in the class, and transcripts of interviews where the instructor describes the instruction that they expect to result in an identified branch of their expected student tree. Instances of routine will be quantified. By triangulating instructor descriptions in interview transcripts, expected student proof constructions, classroom field notes, and expected student discursive maps, I will suggest the portion of instruction that the instructor likely expects to be the source of their students' use of specific routines. The quantity and source of each routine will be compared to student's quantity and sources of routines for each task situation to suggest commonalities and differences in routine selection and the source of routines between students and the instructor in the inter-group analyses of the data. |
|---|---|---|---|

| How do the discourses of students differ from the discourses of their instructors? | How do instructor's expectations of students' injection, surjection, function, and proof discourses prior to and at the end of the course align with students' discourses on injections, surjections, function, and proof as evidenced by their initial and final discursive maps? | Initial discursive maps of students | Mathematical objects and propositions common to the instructor's expected initial and final student discursive maps and mathematical objects and propositions common to students' initial and final discursive maps will be compared using the codes and categories created from the prior constant comparison approach to the analysis of the instructor and student data independently. |
| --- | --- | --- | --- |
| | | Final discursive maps of students | |
| | | Instructor's expected student initial discursive maps | |
| | | Instructor's expected student final discursive maps | |
| | | Video recordings of all instructor interviews | |
| | | Transcripts of all instructor interviews | |
| | For each task situation, how do instructor's expectations of the discourses that students will draw upon to interpret proof task situations involving injective and surjective functions and engage in proof construction | Instructor's expected student initial discursive maps | Mathematical objects and propositions common to the instructor's expected discursive maps for students and expected student trees for each task will be compared to the mathematical objects and propositions common to the student trees for each task as well as students' discursive maps. |
| | | Instructor's expected student final discursive maps | |
| | | Instructor's proof constructions for Tasks A, B, and C | |
| | align with the discourses that students' draw upon as evidenced by their evoked trees and discursive maps? | Instructor's expected student trees for Tasks A, B, and C | |
| | | Present student discursive maps | |
| | | Trees created by students | |
| | | Classroom field notes, notes taken on instructional resources, and alternative resources used by students | |
| | | Video recordings of all interviews | |
| | | Transcripts of all interviews | |
| | How do instructor's expectations of students' function, injection, surjection, and proof discourses based on their reading of students' proof constructions align with students' discourses on function, injections, surjections, and proof as evidenced by their discursive maps? | Classroom field notes, notes taken on instructional resources, and alternative resources used by students | Mathematical objects and propositions that the instructor identifies as likely present or missing within a student's discourse based on their reading of the student's proof construction are compared to the actual student discourse map that was in use at the time of the proof construction. |
| | | Instructors' expected student realization trees and student discursive maps | |
| | | Instructors' narratives written on student proof constructions | |
| | | Students' present discursive maps | |
| | | Students' trees for each task | |

| | | Video recordings of all interviews | |
| --- | --- | --- | --- |
| | | Transcripts of all interviews | |
| | How do instructor's narratives about students' proof and function discourses after reading students' proof constructions align with their narratives on students' proof and function discourses based on students' discursive maps and evoked trees? | Classroom field notes, notes taken on instructional resources, and alternative resources used by students | Prior to the course, I will discuss the narratives that the instructor would like to use to assess the quality of a students' proof and function discourse. As a default example, the instructor can use the narratives: Exceeds Expectations, Meets Expectations, Developing, or Not Meeting Expectations to comment on how they view a students' function and proof discourses when compared to the instructors' course goals. Following each of the student interviews that are conducted during the course, the instructor will examine a series of proof constructions that include de-identified constructions from their students. The instructor will mark each construction with a word or phrase from the agreed upon narratives. As an example, they may identify one construction as "Meets Expectations" and another as "Not Meeting Expectations".<br><br>After the course ends, the instructor is given de-identified student discourse maps labeled with the time frame during the semester that this map was in use by the maps' author. Based solely on the students' function and proof discourse as expressed within the map, the course instructor is asked to mark each map with a word or phrase from the agreed upon narratives.<br><br>After all discursive maps have been marked, the professor is shown one of their students' proof constructions along with the instructors' written narratives. The instructor is then provided with this same student's discursive map that was in use at the time of the proof construction along with the |
| | | Instructors' narratives written on student proof constructions | |
| | | Instructors' narratives written on student discursive concept maps | |
| | | Instructors' expected student trees for each task | |
| | | Students' trees for each task | |
| | | Video recordings of all interviews | |
| | | Transcripts of all interviews | |
| | | | instructor's written narratives. I will produce a table representing the instructors' narratives about student discourses based on student proof constructions and student discursive maps. I will quantify the instructor's consistent vs. inconsistent narratives between artifacts based on differences in ratings provided by the instructor. I will produce a second table representing the instructors' narratives between the three artifacts and quantify consistent vs. inconsistent narratives between artifacts. |

71 Sense-Making in Proof Courses## References

Adiredja, A. P. (2021). Students' struggles with temporal order in the limit definition: uncovering resources using knowledge in pieces. *International Journal of Mathematical Education in Science and Technology*, *52*(9), 1295-1321. doi: 10.1080/0020739X.2020.1754477

Balacheff, N. (1988). A study of students' proving processes at the junior high school level. In I. Wirszup & R. Streit (Eds.), *Proceedings of the 2nd UCSMP International Conference on Mathematics Education*. Reston, VA: National Council of Teachers of Mathematics.

Bansilal, S., Brijlall, O., & Trigueros, M. (2017). An APOS study on pre-service teachers' understanding of injections and surjections. *Journal of Mathematical Behavior, 48*, 22-37.

Barwell, R. (2007). The discursive construction of mathematical thinking: the role of researchers' descriptions. In J. H. Woo, H.C. Lew, K. S. Park, & D. Y. Seo. (Eds.), *Proceedings of the 31st Conference of the International Group for the Psychology of Mathematics Education, Vol. 2* (pp. 49-56). Seoul: PME.

Beckmann, S. (2013). *Mathematics for elementary teachers with activity manual*. New York, NY: Pearson.

Bell, A. (1976). A study of pupils' proof-explanations in mathematical situations. *Educational Studies in Mathematics, 7*(1/2), 23-40.

Breidenbach, D., Dubinsky, E., Hawks, J., & Nichols, D. (1992). Development of the process conception of function. *Educational Studies in Mathematics, 23*, 247-285.

Brown, S. (2017). Who's there? A study of students' reasoning about a proof of existence. *International Journal of Research in Undergraduate Mathematics Education, 3,* 466-495.

74 Sense-Making in Proof Courses

75 Sense-Making in Proof Courses

## Appendix A: Student Clinical Interview Protocol

Good afternoon. My name is Bolanle Salaam and I am a doctoral student in the mathematics education program at the University of Georgia. I am conducting research on student thinking. An important part of this research is to document student discourse development of function and proof throughout their introductory proof courses.

Before we begin, I would like to remind you that I will take every precaution to keep the information that you share during this interview anonymous. Any identifying information that you give during the course of this interview will be replaced with pseudonyms. If at any point in the interview you feel uncomfortable with a question or if you do not feel like answering it, you may choose not to answer that question. You may also end the interview at any time. This interview will take approximately one hour, and it will be helpful if I can record the interview to ensure that I accurately capture your responses. Do I have your permission to record?

*Begin recording if participant consents.*

**Discursive Map Discussion**

1. Prompt: You drew this discursive map for function previously. Can you walk me through how you created this map? *If participant mentions objects that they have not written on the map or does not explain a connection, select the appropriate sub-prompt.*
    a. Probe: You mentioned *signifier-realization*. Were you discussing it in the context of function?
        i. Instruction: If yes, where would that go in your map? Can you draw it?
    b. Prompt: You have *connected signifier-realization to signifier-realization, grouped signifier-realization* under *signifier-realization* in your map. Can you explain what led you to make this connection?
    c. Prompt: I see that you have a link between *signifier-realization and signifier-realization*. That forms *proposition*. Is that what you intend to communicate?
        i. Follow up: Can you describe what you intended to communicate here instead? You can read through as many objects and links as you like. I just want to know what you are trying to say about *signifier-realization* here.
    *Repeat for all portions of the map that have not been explained by the participant.*

2. Prompt: Would you like to add anything else or change the format of your map? You can write directly on your map.
    a. Probe: You have *added to, removed, connected signifier to, grouped, split up* your *signifier-realization* and *describe other signifiers and realizations that were associated with this change and the nature of change.* Can you explain what led you to make this change?
        i. Follow up: You view *list each proposition resulting from this change* as a true statement. Am I correct? If not, please correct me.

      b. Follow up: Previously, you *describe state of signifier-realization in previous map* because *repeat uses described by participant in their prior map,* but now you *describe current connections and repeat uses for object.* Did I summarize that accurately?

      c. Prompt: So, *describe change* is now *describe positioning on map.* This places it *describe global map impact of change if participant does not describe it*. Am I correct that this is how you view this?

*Transition:* Let me know if you realize that you want to add something to your map later on in the interview. For now, I'd like to talk about some of your experiences.

**Student Experiences Discussion**

3. Prompt: Can you tell me about when you were first introduced to the term *function?*
   a. Probe: Was this introduction in the school setting?
      i. Probe: If so, what grade or course were you in?
      ii. Probe: If not, can you describe the circumstances around this introduction to *function?*
   b. Follow up: Walk me through other experiences you have had with the term "function".
   c. Prompt: What types of experiences did you include in your *function* map? Is there anything that led you to exclude some of your experiences?
   d. Prompt: Is there anything else you would like me to know about your experiences with *function?*
4. Prompt (Optional based on time and student experience): Can you tell me about when you were first introduced to the term *proof*?
   a. Probe: Was this introduction in the school setting?
      i. Probe: If so, what grade or course were you in?
      ii. Probe: If not, can you describe the circumstances around this introduction to *proof?*
   b. Follow up: Walk me through other experiences you have had with the term *proof*.
   c. Prompt: Is there anything else you would like me to know about your experiences with *proof?*

*Final Question:* Excellent. Is there anything else that you would like me to know about your map or anything you want to change? We can also take this time to make a clean final map.

*Wrap-Up:* Thank you for discussing your thinking and past experiences with me! I will contact you if you are selected for further interviews. Do not forget to complete the follow-up survey within 24 hours to get your *incentive.* Your responses will not impact you receiving your *incentive* so please be as honest as possible.

## Appendix B: Instructor Clinical Interview Protocol

Good afternoon. My name is Bolanle Salaam and I am a doctoral student in the mathematics education program at the University of Georgia. I am conducting research on student thinking and discourse. An important part of this research is to document student discourse development of proof and function before, during, and at the end of their introductory proof courses. It would be helpful for me to know how you think about function, proof, as well as your thoughts on student thinking in this course.

Before we begin, I would like to remind you that I will take every precaution to keep the information that you share during this interview anonymous. Any identifying information that you give during the course of this interview will be replaced with pseudonyms. If at any point in the interview you feel uncomfortable with a question or if you do not feel like answering it, you may choose not to answer that question. You may also end the interview at any time. This interview will take approximately one hour, and it will be helpful if I can record the interview to ensure that I accurately capture your responses. Do I have your permission to record?

*Begin recording if participant consents.*

**Discursive Map Discussion**

1. Prompt: You drew this discursive map for function previously. Can you walk me through how you created this map? *If participant mentions objects that they have not written on the map or does not explain a connection, select the appropriate sub-prompt.*
    a. Probe: You mentioned *signifier-realization*. Were you discussing it in the context of function?
        i. Instruction: If yes, where would that go in your map? Can you draw it?
    b. Prompt: You have *connected signifier-realization to signifier-realization, grouped signifier-realization* under *signifier-realization* in your map. Can you explain what led you to make this connection?
    c. Prompt: I see that you have a link between *signifier-realization and signifier-realization.* That forms *proposition.* Is that what you intend to communicate?
        i. Follow up: Can you describe what you intended to communicate here instead? You can read through as many objects and links as you like. I just want to know what you are trying to say about *signifier-realization* here.
    *Repeat for all portions of the map that have not been explained by the participant.*

2. Prompt: Would you like to add anything else or change the format of your map? You can write directly on your map.
    a. Probe: You have *added to, removed, connected signifier to, grouped, split up* your *signifier-realization* and *describe other signifiers and realizations that were associated with this change and the nature of change.* Can you explain what led you to make this change?

        ii. Follow up: So, you view *list each proposition resulting from this change* as a true statement. Am I correct? If not, please correct me.
    b. Follow up: Previously, you *describe state of signifier-realization in previous map* because *repeat uses described by participant in their prior map,* but now you *describe current connections and repeat uses for object.* Did I summarize that accurately?
    c. Prompt: So, *describe change* is now *describe positioning on map.* This places it *describe global map impact of change if participant does not describe it.* Am I correct that this is how you view this?

3. Prompt: You were asked to draw a discursive map that you feel reflects how you expect your typical incoming students think of function as well as an expected discursive map that reflects your goals for how students' function discourse will be organized by the time they finish this course. These are the expected discourse maps of students that you drew. Let's begin with the initial map. Can you walk me through how you created this map?

    a. Probe: In which course do you feel that students would have encountered *signifier-realization or map portion* in your map?
        i. Follow up: Is this when you first encountered *signifier-realization or map portion* as well? If you remember, when did you first encounter *signifier-realization or map portion*?
        ii. Follow up: What are your experiences with student preparation in *signifier-realization or map portion* when they enter this course?

    *Begin with function and then ask about other signifier-realizations present in the instructor's expectation map until all propositions are discussed.*

4. Prompt: If possible, can you describe how you came to the map of *function* that you have now? At what point do you believe you achieved this level of organization? What about the expected final discursive map of students that you drew for *function*?

*Transition:* Now, I would like to understand your thoughts about the introductory proof course that you will be teaching.

**Introductory Proof Course Discussion (Alternative: Distribute these questions in a survey.)**

5. In your own words, why should students take this course?
    a. Follow up: Do you currently feel that this course accomplishes the aims that you have discussed?
6. How do you expect student's proof discourse to be organized by the end of the semester? You may draw a discursive map for proof if you wish. *The remaining questions can be* *skipped if they are addressed in the instructor's response to question 5 or 6. The instructor is not required to draw a proof map if they do not feel it is necessary,*
7. Which topics do you feel students have the most success with in this course?
    a. Follow up: What are you looking for students to do in regard to *topic* to consider them successful?
8. Which topics do you feel pose the most difficulty for students in this course?
    a. Follow up: What do you feel is the source of student difficulty with *topic?*
9. What is your experience with students in this course regarding function, injections, and surjections?
10. How does your course address student difficulty with proof construction?

*Final Question:* Is there anything else that you would like me to know about your course or transition to proof courses in general?

*Wrap-Up:* Thank you for discussing your thinking and past experiences with me! After I have students complete a task, I will send the task over to you. As a reminder, you will write at least one proof construction that reflects how you expect your students to approach the task, and you will also draw a tree that reflects how you would expect your students to parse this task before ultimately arriving at their construction. If you anticipate other common approaches, you are free to write additional constructions and trees if you like. I will collect these from you in advance of my next interview with you. I look forward to reading your responses to the student tasks.

## Appendix C: Student Task-Based Interview Protocol (Task A)

Welcome back and thank you for participating! As a reminder, you may end your participation at any time. You may also decline to answer any questions that you are not comfortable with answering. It would be beneficial for me to record this interview. May I record?

*Record if participant consents.*

**Map Changes**

I am going to show you the map you approved after our last meeting. Let me know if you want to make any changes. *Show participant their function map.*

1. Prompt: I see that you *describe map changes generally*. Did I miss any of your changes? *If no changes are made to student map(s), continue to proof and proving discussion.*

    a. Prompt: You have *added to, removed, connected signifier to, grouped, split up* your *signifier-realization* and *describe other signifiers and realizations that were associated with this change and the nature of change.* Can you explain what led you to make this change?
        i. Follow up: So, you view *list each proposition resulting from this change* as a true statement. Am I correct? If not, please correct me.
    b. Follow up: Previously, you *describe state of signifier-realization in previous map* because *repeat uses described by participant in their prior map,* but now you *describe current connections and repeat described uses for object.* Did I summarize that accurately?
    c. Prompt: That's interesting. Can you remember where you *repeat student use of object* or how you came to know that *describe new proposition made by changes*?
        i. Follow up for non-instructional materials: You stated that you used *non-instructional resource.* Do you mind sharing that resource with me?
        ii. Follow up: What led you to *non-instructional resource* as opposed to *name instructional materials that address the content of the change*?
        iii. Follow up: Your summary of *non-instructional resource* is that *describe new proposition* because *describe features of non-instructional resource that participant cites as the source of the map change*? Am I correct? If not, please correct me.
    d. Prompt: So, *describe signifier-realization that has changed* is now *describe positioning on map.* This places it *describe global map impact of change if participant does not describe it*. Am I correct that this is how you view this?
    e. Prompt: Is there anything else that you would like to tell me about your changes? Is there anything that we have missed?

    *Repeat until all changes have been discussed by participant. If participant makes additional changes, repeat sub-protocol for the new change.*

**Proof Discursive Map**

Now, I would like you to make a discursive map for proof. You have learned about these for about a month now. Take about 10 minutes to write as much as you can about how you think about proof. Organize it into a discursive map. Here is your function map and some materials. Let me know when you are done. *Provide participant with their function discursive map, blank sheets of paper, and writing materials. Do not indicate that proof map must be drawn on a new sheet nor connected to the function map. If participant asks where the map should be drawn, state that it is up to them and where they feel comfortable drawing it as long as it reflects how they think about proof.*

1. Prompt: Great! As a note, you do not have to *draw separate maps, combine both maps, link both maps* you also could have *draw separate maps, combine both maps, link both maps* if that is how you think about them. Is your current map(s) an accurate depiction of how you think about *function* and *proof*? *If participant says no, provide instruction.*
    a. Instruction: Feel free to make changes as you see fit.

2. Prompt: Now that you have made your map(s), can you walk me through how you created this map? *If participant mentions objects that they have not written on the map or does not explain a connection, select the appropriate sub-prompt.*
    a. Probe: You mentioned *signifier-realization*. Were you discussing it in the context of *proof, function, proof and function*?
        i. Instruction: If yes, where would that go in your map? Can you draw it?
    b. Prompt: You have *connected signifier-realization to signifier-realization, grouped signifier-realization* under *signifier-realization* in your map. Can you explain what led you to make this connection?
    c. Prompt: I see that you have a link between *signifier-realization and signifier-realization.* That forms *proposition.* Is that what you intend to communicate?
        i. Follow up: Can you describe what you intended to communicate here instead? You can read through as many objects and links as you like. I just want to know what you are trying to say about *signifier-realization* here.
    d. Prompt: Can you remember where you *repeat student use of object* or how you came to *describe proposition in map*?
        i. Follow up for non-instructional materials: You stated that you used *non-instructional resource.* Do you mind sharing that resource with me?
        ii. Follow up: What led you to *non-instructional resource* as opposed to *name instructional materials that address the content of the change*?
        iii. Follow up: So, your summary of *non-instructional resource* is that *describe new proposition* because *describe features of non-instructional resource that participant cites as the source of the map change*? Am I correct? If not, please correct me.

*Repeat for all portions of the map that have not been explained by the participant.*

3. Prompt: Would you like to add anything else or change the format of your map? You can write directly on your map.
   a. Probe: You have *added to, removed, connected signifier to, grouped, split up* your *signifier-realization* and *describe other signifiers and realizations that were associated with this change and the nature of change.* Can you explain what led you to make this change?
      i. Follow up: So, you view *list each proposition resulting from this change* as a true statement. Am I correct? If not, please correct me.
   b. Follow up: So previously, you *describe state of signifier-realization in previous map* because *repeat uses described by participant in their prior map,* but now you *describe current connections and repeat uses for object.* Did I summarize that accurately?
   c. Prompt: So, *describe change* is now *describe positioning on map.* This places it *describe global map impact of change if participant does not describe it.* Am I correct that this is how you view this?

*Transition:* Let me know if you realize that you want to add something to your map later on in the interview.

**Task-Situation (Sample Task from Task Bank)**

1. Ok. Now I am going to have you work on a task. Try your best to communicate everything you are thinking. It is really important that you let me know what you see in this task that is leading you to act. If you are feeling a bit unsure about where to start or if you need to spend more time thinking about this, we can start building your tree first to help you explore what this task is asking you to do and possible ways that you can respond. You can do your proof construction afterward. *Give student task situation. Do not provide definitions for this task. Do not allow textbook use.*

**REDACTED IMAGE**

**REDACTED CITATION**

Alternate Sample Task:

*Prove that for* $f: \mathbb{R} \to \mathbb{R}, f(s) = \dfrac{1}{2^s}$ *satisfies* $f(x) = f(y) \Rightarrow x = y$ *for all* $x, y$ *in the domain.*

Acceptable probes for evoking student discourse:

- Try re-reading the task. Is there any portion of the task that you are focusing on?
- What are you thinking about currently? Where do you feel stuck?
- When you get stuck in this class on assignments, what do you usually do?
- When you have more than one option that you can explore in a problem, how do you usually choose from your options in this class?
    - Why don't we list the options you are currently considering by writing them in your tree? That may help you process the task.

This probe should only be used if no progress has been made by the final twenty minutes of this interview and the acceptable probes have been exhausted:

- You have said that *repeat any if-then statements participant has made*. Do you feel that *repeat "if" portion* is the only way to do this problem or are there other possibilities that you are considering? *Choose the follow-up question to ask first based on participant's expressed inclinations.*
    - Follow up: Let's say you did know *repeat the "if" portion of participant's statement*, how would you complete the task at that point?
    - Follow up: Now let's say you know that *false version of prior hypothetical statement*, what would you do at that point?

*Transition: (a) If student does not complete task, proceed to wrap-up. (b) Otherwise, continue.*

2. Prompt: Ok. Can you tell me what parts of this problem made you think of your approach? *Make note of participant's expressed task prompts, routines, and sub-routines. Ask about similar simple tasks to identify regularities in student's routine execution or inhibition. Allow participant to create tree while describing their approach if this has not been done already.*

    a. Probe: Which parts of the task did you notice first? *Make note of the participant's explication of possible key words used by this participant for a routine.*
        i. Follow up: What did that part of the task make you think of?
    b. Probe: Is this how you always interpret this *word/symbol/compound phrase*? *Make note of the participant's recollection of additional evoked uses for key words and any applicability conditions for their expressed interpretation.*
        i. Follow up: What about in this course specifically? Is there some other scenario where you would have a different interpretation?
        ii. (Optional) Follow up: Do the *keywords described by participant* have to appear together for you to do this routine? Are there other features of similar tasks that you need to examine before you would know that you can *describe routine*?
    c. Probe: What are some other types of problems that would make you proceed in the way you did for this problem, if any? *Make note of key words and applicability conditions that participant states lead them to enact the same routine used in this task.*

      d. Prompt: What criteria did you use to decide that you were done with this proof construction? *Make note of the closing conditions identified by the participant.*
      e. Prompt: You have stated that you saw *key words* and that they made you *describe routine identified by the student* and that you knew you were done when *describe closing conditions identified by the participant.* Is this correct?
      f. Prompt: Are there any portions of the task-situation that you find irrelevant to how you approached this problem?
          i. Follow up: I see *list unused objects in task-situation* in the task situation. Were any of these relevant to how you approached this problem? If so, please explain.
      g. Prompt: Here is *element of construction* that I see in your proof construction. Can you let me know what led you to write this?
          i. Follow up: I don't see where you discuss *element of construction* in your tree. Can you explain how you arrived here some more? We will update your tree to reflect this.

*Repeat b-g for any sub-routines the participant describes.*

3. Prompt: Ok. Here is your discursive map(s). Are there any elements of your map(s) that you used to do this task? If so, show me. *If participant can trace their complete approach using the map, have them draw a tree that copies the relevant parts of their discursive map if they have not created a tree already; they should include any decision-making discourse they have used. If the participant's routine is not explicated by the discursive map, have student create tree at this point if their tree has not been created already.*

    a. Prompt: I see that you have a link between *mathematical object* and *mathematical object.* That forms *proposition.* Is that what you intend to communicate? *Ask about any propositions that appear to contradict or support the participant's proof construction and tree approach.*
      i. Follow up: Can you describe what you intend to communicate here? You can read through as many objects and links as you like. I just want to know what you are trying to say about *mathematical object* here.

*Final Question:* Did you have any difficulty answering any of my questions about how you constructed your proof? Do you have comments on how you approached the proof? We can also take this time to make clean versions of your final map and tree.

*Wrap-Up:* Thank you for discussing your thinking with me! I will contact you if you are selected for further interviews. Do not forget to complete the follow-up survey within 24 hours to get your *incentive.* Your responses will not impact you receiving your *incentive* so please be as honest as possible.

## Appendix D: Student Task-Based Interview Protocol (Task B)

Welcome back and thank you for participating! As a reminder, you may end your participation at any time. You may also decline to answer any questions that you are not comfortable with answering. It will be beneficial for me to record this interview. May I record?

*Record if participant consents.*

**Map Changes**

I am going to show you the map you approved after our last meeting. Let me know if you want to make any changes. *Show participant their map(s).*

1. Prompt: I see that you *describe each map change.* Did I miss any of your changes? *If no changes are made to student map(s), continue to task-situation.*

    a. Prompt: You have *added to, removed, connected signifier to, grouped, split up* your *signifier-realization* and *describe other signifiers and realizations that were associated with this change and the nature of change.* Can you explain what led you to make this change?
        i. Follow up: So, you view *list each proposition resulting from this change* as a true statement. Am I correct? If not, please correct me.
    b. Follow up: Previously, you *describe state of signifier-realization in previous map* because *repeat uses described by participant in their prior map,* but now you *describe current connections and repeat uses for object.* Did I summarize that accurately?
    c. Prompt: That's interesting. Can you remember where you *repeat student use of signifier-realization* or how you came to *describe new proposition made by changes*?
        i. Follow up for non-instructional materials: You stated that you used *non-instructional resource.* Do you mind sharing that resource with me?
        ii. Follow up: What led you to *non-instructional resource* as opposed to *name instructional materials that address the content of the change*?
        iii. Follow up: So, your summary of *non-instructional resource* is that *describe new proposition* because *describe features of non-instructional resource that participant cites as the source of the map change*? Am I correct? If not, please correct me.
    d. Prompt: So, *describe change* is now *describe positioning on map.* This places it *describe global map impact of change if participant does not describe it*. Am I correct that this is how you view this?
    e. Prompt: Is there anything else that you would like to tell me about your changes? Is there anything that we have missed?

    *Repeat until all changes have been discussed by participant. If participant makes additional changes, repeat sub-protocol for the new change.*

**Task-Situation (Actual Task)**

1. Ok. Now I am going to have you work on a task. Try your best to communicate everything you are thinking. It is really important that you let me know what you see in this task that is leading you to act. If you are feeling a bit unsure about where to start or if you need to spend more time thinking about this, we can start building your tree first to help you explore what this task is asking you to do and possible ways that you can respond. You can do your proof construction afterward. *Give student task situation. Do not provide definitions for this task unless asked for by student. If student asks for definitions, note the time on task before this occurred. If student spends ten minutes on this task without making initial progress, provide definitions at that time. Content-specific help outside of definitions should not be provided. The student may appeal to their course textbook and course materials if they ask for them, but these must be documented in detail by interviewer.*

$$\text{Prove or disprove that } \frac{1}{2^x} \text{ is injective and/or surjective from } \mathbb{R} \to \mathbb{R}.$$

Acceptable probes for evoking student discourse:
- Try re-reading the task. Is there any portion of the task that you are focusing on?
- What are you thinking about currently? Where do you feel stuck?
- When you get stuck in this class on assignments, what do you usually do?
- When you have more than one option that you can explore in a problem, how do you usually choose from your options in this class?
    - Why don't we list the options you are currently considering by writing them in your tree? That may help you process the task.

This probe should only be used if no progress has been made by the final twenty minutes of this interview and the acceptable probes have been exhausted:
- You have said that *repeat any if-then statements participant has made*. Do you feel that *repeat "if" portion* is the only way to do this problem or are there other possibilities that you are considering? *Choose the follow-up question to ask first based on participant's expressed inclinations.*
    - Follow up: Let's say you did know *repeat the "if" portion of participant's statement*, how would you complete the task at that point?
    - Follow up: Now let's say you know that *false version of hypothetical statement*, what would you do at that point?

*Transition: (a) If student does not complete task, proceed to wrap-up. (b) Otherwise, continue.*

2. Prompt: Ok. Can you tell me what parts of this problem made you think of your approach? *Make note of participant's expressed task prompts, routines, and sub-routines. Ask about similar simple tasks to identify regularities in student's routine* *execution or inhibition. Allow participant to create tree while describing their approach if this has not been done already.*

   a. Probe: Which parts of the task did you notice first? *Make note of the participant's explication of possible key words used by this participant for a routine.*
      i. Follow up: What did that part of the task make you think of?
   b. Probe: Is this how you always interpret this *word/symbol/compound phrase*? *Make note of the participant's recollection of additional evoked uses for key words and any applicability conditions for their expressed interpretation.*
      i. Follow up: What about in this course specifically? Is there some other scenario where you would have a different interpretation?
      ii. (Optional) Follow up: Do the *keywords described by participant* have to appear together for you to do this routine? Are there other features of similar tasks that you need to examine before you would know that you can *describe routine*?
   c. Probe: What are some other types of problems that would make you proceed in the way you did for this problem, if any? *Make note of key words and applicability conditions that participant states lead them to enact the same routine used in this task.*
   d. Prompt: What criteria did you use to decide that you were done with this proof construction? *Make note of the closing conditions identified by the participant.*
   e. Prompt: You have stated that you saw *key words* and that they made you *describe routine identified by the student* and that you knew you were done when *describe closing conditions identified by the participant.* Is this correct?
   f. Prompt: Are there any portions of the task-situation that you find irrelevant to how you approached this problem?
      i. Follow up: I see *list unused objects in task-situation* in the task situation. Were any of these relevant to how you approached this problem? If so, please explain.
   g. Prompt: Here is *element of construction* that I see in your proof construction. Can you let me know what led you to write this?
      i. Follow up: I don't see where you discuss *element of construction* in your tree. Can you explain how you arrived here some more? We will update your tree to reflect this.

   *Repeat b-g for any sub-routines the participant describes.*

3. Prompt: Ok. Here is your discursive map(s). Are there any elements of your map(s) that you used to do this task? If so, show me. *If participant can trace their complete approach using the map, have them draw a tree that copies the relevant parts of their discursive map if they have not created a tree already; they should include any decision-making* *discourse they have used. If the participant's routine is not explicated by the discursive map, have student create tree at this point if their tree has not been created already.*

   a. Prompt: I see that you have a link between *mathematical object* and *mathematical object*. That forms *proposition*. Is that what you intend to communicate? *Ask about any propositions that appear to contradict or support the participant's proof construction and tree approach.*
        i. Follow up: Can you describe what you intend to communicate here? You can read through as many objects and links as you like. I just want to know what you are trying to say about *mathematical object* here.

*Final Question:* Did you have any difficulty answering any of my questions about how you constructed your proof? Do you have comments on how you approached the proof? We can also take this time to make a clean final map and tree.

*Wrap-Up:* Thank you for discussing your thinking with me! I will contact you if you are selected for further interviews. Do not forget to complete the follow-up survey within 24 hours to get your *incentive.* Your responses will not impact you receiving your *incentive* so please be as honest as possible.

**Appendix E: Student Task-Based Interview Protocol (Task C)**

Welcome back and thank you for participating! This will be the final student interview. As a reminder, you may end your participation at any time. You may also decline to answer any questions that you are not comfortable with answering. It will be beneficial for me to record this interview. May I record?

*Record if participant consents.*

**Map Changes**

I am going to show you the map you approved after our last meeting. Let me know if you want to make any changes. *Show participant their map(s).*

1. Prompt: I see that you *describe each map change.* Did I miss any of your changes? *If no changes are made to student map(s), continue to task-situation.*

    a. Prompt: You have *added to, removed, connected signifier to, grouped, split up* your *signifier-realization* and *describe other signifiers and realizations that were associated with this change and the nature of change.* Can you explain what led you to make this change?
        i. Follow up: So, you view *list each proposition resulting from this change* as a true statement. Am I correct? If not, please correct me.
    b. Follow up: So previously, you *describe state of signifier-realization in previous map* because *repeat participant's described uses in their prior map,* but now you *describe current connections and repeat uses for object.* Did I summarize that accurately?
    c. Prompt: That's interesting. Can you remember where you *repeat student use of object* or how you came to *describe new proposition made by changes*?
        i. Follow up for non-instructional materials: You stated that you used *non-instructional resource.* Do you mind sharing that resource with me?
        ii. Follow up: What led you to *non-instructional resource* as opposed to *name instructional materials that address the content of the change*?
        iii. Follow up: So, your summary of *non-instructional resource* is that *describe new proposition* because *describe features of non-instructional resource that participant cites as the source of the map change*? Am I correct? If not, please correct me.
    d. Prompt: So, *describe change* is now *describe positioning on map.* This places it *describe global map impact of change if participant does not describe it*. Am I correct that this is how you view this?
    e. Prompt: Is there anything else that you would like to tell me about your changes? Is there anything that we have missed?

    *Repeat until all changes have been discussed by participant. If participant makes additional changes, repeat sub-protocol for the new change.*

**Task-Situation** (Samples from Task Bank)

1. Ok. Now I am going to have you work on a task. Try your best to communicate everything you are thinking. It is really important that you let me know what you see in this task that is leading you to act. If you are feeling a bit unsure about where to start or if you need to spend more time thinking about this, we can start building your tree first to help you explore what this task is asking you to do and possible ways that you can respond. You can do your proof construction afterward. *Give student task situation. Do not provide definitions for this task unless asked for by student. If student asks for definitions, note the time on task before this occurred. If student spends ten minutes on this task without making initial progress, provide definitions at that time. Content-specific help outside of definitions should not be provided. The student may appeal to their course textbook and materials if they ask for them, but these must be documented in detail by interviewer.*

*Alternate Sample Task*: *Prove or disprove that $h(n) = 3n$ is injective and/or surjective from $\mathbb{Z} \to \mathbb{Z}$.*
(Thoma & Nardi, 2018, p. 169)

*Alternate Sample Task*: *Let $f, g$ be two functions whose domains and ranges are subsets of the set of real numbers. Prove or find a counterexample to the following statements:*
*(a) If $f, g$ are both injective, then it follows that $f + g$ is injective.*
*(b) If $f, g$ are both injective, then it follows that $f + g$ is surjective.*
(Breidenbach et. al., 1992, p. 283)

*Alternate Sample Task*: *Problems with variations on maps between $\mathbb{R}, \mathbb{N}, \mathbb{Z} \to \mathbb{R}, \mathbb{N}, \mathbb{Z}$*

Acceptable probes for evoking student discourse:
- Try re-reading the task. Is there any portion of the task that you are focusing on?
- What are you thinking about currently? Where do you feel stuck?
- When you get stuck in this class on assignments, what do you usually do?
- When you have more than one option that you can explore in a problem, how do you usually choose from your options in this class?
    - Why don't we list the options you are currently considering by writing them in your tree? That may help you process the task.

This probe should only be used if no progress has been made by the final twenty minutes of this interview and the acceptable probes have been exhausted:

- You have said that *repeat any if-then statements participant has made*. Do you feel that *repeat "if" portion* is the only way to do this problem or are there other possibilities that you are considering? *Choose the follow-up question to ask first based on participant's expressed inclinations.*
    - Follow up: Let's say you did know *repeat the "if" portion of participant's statement*, how would you complete the task at that point?
    - Follow up: Now let's say you know that *false version of hypothetical statement*, what would you do at that point?

*Transition: (a) If student does not complete task, proceed to wrap-up. (b) Otherwise, continue.*

2. Prompt: Ok. Can you tell me what parts of this problem made you think of your approach? *Make note of participant's expressed task prompts, routines, and sub-routines. Ask about similar simple tasks to identify regularities in student's routine execution or inhibition. Allow participant to create tree while describing their approach if this has not been done already.*

    a. Probe: Which parts of the task did you notice first? *Make note of the participant's explication of possible key words used by this participant for a routine.*
        i. Follow up: What did that part of the task make you think of?
    b. Probe: Is this how you always interpret this *word/symbol/compound phrase*? *Make note of the participant's recollection of additional evoked uses for key words and any applicability conditions for their expressed interpretation.*
        i. Follow up: What about in this course specifically? Is there some other scenario where you would have a different interpretation?
        ii. (Optional) Follow up: Do the *keywords described by participant* have to appear together for you to do this routine? Are there other features of similar tasks that you need to examine before you would know that you can *describe routine*?
    c. Probe: What are some other types of problems that would make you proceed in the way you did for this problem, if any? *Make note of key words and applicability conditions that participant states lead them to enact the same routine used in this task.*
    d. Prompt: What criteria did you use to decide that you were done with this proof construction? *Make note of the closing conditions identified by the participant.*
    e. Prompt: You have stated that you saw *key words* and that they made you *describe routine identified by the student* and that you knew you were done when *describe closing conditions identified by the participant.* Is this correct?
    f. Prompt: Are there any portions of the task-situation that you find irrelevant to how you approached this problem?

    i. Follow up: I see *list unused objects in task-situation* in the task situation. Were any of these relevant to how you approached this problem? If so, please explain.
  g. Prompt: Here is *element of construction* that I see in your proof construction. Can you let me know what led you to write this?
    i. Follow up: I don't see where you discuss *element of construction* in your tree. Can you explain how you arrived here some more? We will update your tree to reflect this.

*Repeat b-g for any sub-routines the participant describes.*

3. Prompt: Ok. Here is your discursive map(s). Are there any elements of your map(s) that you used to do this task? If so, show me. *If participant can trace their complete approach using the map, have them draw a tree that copies the relevant parts of their discursive map if they have not created a tree already; they should include any decision-making discourse they have used. If the participant's routine is not explicated by the discursive map, have student create tree at this point if their tree has not been created already.*

  a. Prompt: I see that you have a link between *mathematical object* and *mathematical object*. That forms *proposition*. Is that what you intend to communicate? *Ask about any propositions that appear to contradict or support the participant's proof construction and tree approach.*
    i. Follow up: Can you describe what you intend to communicate here? You can read through as many objects and links as you like. I just want to know what you are trying to say about *mathematical object* here.

*Final Question:* Did you have any difficulty answering any of my questions about how you constructed your proof? Do you have comments on how you approached the proof? We can also take this time to make a clean final map and tree.

*Wrap-Up:* Thank you for discussing your thinking with me and for participating in these interviews over the semester. Do not forget to complete the follow-up survey within 24 hours to get your final *incentive*. Your responses will not impact you receiving your *incentive* so please be as honest as possible.

# Appendix F: Instructor Task Based Interviews Protocol

*This protocol will be used for all three task-based interviews that occur during the semester in response to Task A, Task B, and Task C.*

Before we begin, I would like to remind you that I will take every precaution to keep the information that you share during this interview confidential. Any identifying information that you give during the course of this interview will be replaced with pseudonyms. If at any point in the interview you feel uncomfortable with a question or if you do not feel like answering it, you may choose not to answer that question. You may also end the interview at any time. This interview will take approximately one hour, and it will be helpful if I can record the interview to ensure that I accurately capture your responses. Do I have your permission to record?

*Begin recording if participant consents.*

**Expected Construction and Trees Discussion**

The purpose of today is for you to explain your proof construction and the trees that you thought would be most like the student approaches to *specify task.* You will then be given the opportunity to examine de-identified sample constructions. Some of these constructions will be from your classroom, while other constructions were contributed by participants that are not in your course, graduate students, mathematics teachers, and mathematics instructors. Contributors were asked to write their constructions at the level they would expect of a student in an introductory proof course. You will not be informed of who made the sample constructions that you will view.

1. Prompt: Can you walk me through your proof construction for *specify task*?
    a. Prompt: Does your expected construction represent what you would deem as an acceptable construction for this task? *If yes, continue to 1b.*
        i. Follow up: If not, please describe the reasoning that went into the construction you expect of your students on this task.
    b. Prompt: Based on the instruction that was completed by the time of interviews, do you believe your students had received the instruction necessary to provide a construction to this task that you would deem acceptable?
    c. Prompt: Based on the instruction that was completed by the time of interviews, what aspects of this task do you think students would have struggled with, if any?
    d. Prompt: Based on the instruction that was completed by the time of interviews, what aspects of this task do you think should have presented minimal difficulty for your students?
    e. Follow up: Do you have any comments on your expected proof construction or the conversation we have had in this interview so far?

2. Prompt: Can you walk me through your expected student trees for this task? Here are copies of the expected initial and final student discursive maps that you drew at the beginning of the semester. Feel free to use them as necessary as you explain your reasoning.

    a. Probe: Which parts of the task do you expect students to notice first? *Make note of the participant's explication of possible key words used by students for a routine.*
        i. Follow up: What did you expect that portion of the task to make students think of?
    b. Probe: Is this how you expect students to interpret this *word/symbol/compound phrase* in all mathematical situations? *Make note of the participant's recollection of additional evoked uses for key words.*
        i. Follow up: What about in this course specifically? Is there some other scenario where you would expect students to interpret *word, symbol, or image* differently?
        ii. (Optional. Use insider judgment on phrasing.) Follow up: Do the *keywords described by participant* have to appear together for you to expect students to do this routine? Are there other features of similar tasks that you would expect students to examine before they enact the same routine?
    c. Probe: Are there other types of problems/families of problems that you would expect students to complete by applying the same routine that you used in you expected student tree? *Make note of key words and applicability conditions that participant states should lead students to enact the same routine used in this task.*
    d. Prompt: What criteria do you expect students to use to decide that they have finished this proof construction and no further work is needed? *Make note of the closing conditions identified by the participant.*
    e. Prompt: You have stated that you expected students to initially notice *key words* and that this should lead them to *describe expected routine identified by the instructor* and that they would know they were done with the construction by *describe closing conditions identified by the participant.* Is this correct?
    f. Prompt: Are there any portions of the task-situation that were irrelevant to how you expected students to approach this task?
        i. Follow up: I see *list unused objects in task-situation* in the task situation. Were any of these relevant to how you expected students to approach this task? If so, please explain.
    g. Prompt: Here is *element of construction* that I see in your expected proof construction. Can you let me know what led you to write this?
        i. Follow up: I don't see where you discuss *element of construction* in your expected student tree. Can you explain how you expected students to *discuss unexplained aspect of construction* some more? We will update your tree to reflect this.

    *Repeat b-g for any sub-routines the participant describes.*

3. Prompt: Do you have any other comments about this task or how you expect students to approach it?

**Blind Proof Construction Discussion** (Expected Timing: 30 minutes for Three Constructions)

Now I will show you a series of proof constructions that were written in response to this task. As a reminder, you may or may not see proof constructions from your classroom in this sample. I would like for you to read through each construction, comment on the quality of the construction, and specify how you would assess this construction if one of your students were to hand this in on an assignment using our agreed upon rating scale. I would also like for you to comment on how you believe the author of the construction may have thought about the task and how their function and/or proof discourse may be organized. You can use the expected student discursive map and trees that you created to do this if you wish. I will have you write your assessment on each construction.

1. Here is *construction sample* in response to *specify task*. Please take a moment to read the proof construction and comment on the quality of the construction.
    a. Probe: What are you using as the basis for determining the quality of this construction?
    b. Probe: How do you think the author interpreted this task? How do you think this task was interpreted? You can compare this to the expected tree you made for this task if it helps you explain your thoughts.
    c. Probe: How would you expect a student with this construction's function discourse to be organized? Are there any aspects of function discourse that you would not expect to be present in this student's discourse?
    d. Probe: Do you have any other comments about this construction?

    *Ask instructor to write their assessment on the construction. Repeat sub-protocol for all sample constructions.*

*Final Question:* Do you have any comments for me regarding this task, your expected construction, the associated tree, or the sample constructions that we read today?

*Wrap-Up:* Thank you for discussing your thinking with me! I look forward to speaking with you at the next interview.

# Appendix G: Instructor Final Clinical Interview Protocol

Before we begin, I would like to remind you that I will take every precaution to keep the information that you share during this interview confidential. Any identifying information that you give during the course of this interview will be replaced with pseudonyms. If at any point in the interview you feel uncomfortable with a question or if you do not feel like answering it, you may choose not to answer that question. You may also end the interview at any time. This interview will take approximately one hour, and it will be helpful if I can record the interview to ensure that I accurately capture your responses. Do I have your permission to record?

*Begin recording if participant consents.*

**Instructor Discursive Map Discussion**

1. Prompt: You drew this discursive map for your function discourse previously. Do you have any changes to this map? *If participant mentions objects that they have not written on the map, ask them to include that object in how they think about function.*
    a. Probe: Is there anything that led you to make the new connections on your map?
2. Prompt: You drew these expected initial and final student discursive maps for *function* previously. Based on how the course went this semester and your experience with your students, do you want to make any changes to how you think their final student discursive concept maps will be organized?
    a. Probe: What led you to *make new connection, reorganize map, remove connection*?

**Student Discursive Concept Maps Discussion**

During the semester, I showed you a series of proof constructions that students in your course completed as well as contributed maps from participants in a prior phase of the study and volunteers. I have collected the trees and discursive maps for each construction that you examined that was made by one of your students. Only your students are included in this sample.

I will show you representative function discursive maps made by students in this study throughout the semester. If a student made a separate proof discursive map, I will show you that as well. There is no order to these maps, and they do not align with the order that I used for showing you proof constructions throughout the semester. I would like for you to comment on how you think a student with each of these maps would have performed on each task and any approaches you think they would have used to access the task. You can use elements from the sample student's discursive map(s) to clarify your reasoning. You will write your commentary on the maps.

1. Prompt: Here is *one student's map(s) or a representative discursive map sample*. Please take a moment to read the map and comment on the quality of the map, any concerns that you have, and place a rating on the map using our agreed upon scale.
    a. Probe: What are you using as the basis for determining the quality of this map?

      b. Probe: How do you think this student interpreted and performed on *specify task*? How do you think this task was interpreted by students with similar maps?
      c. Do you have any other comments about this map?

*Have instructor label discursive maps. Repeat protocol for multiple discursive maps. Group maps by the task where the similar map was used by the student participant.*

**Student Realization Tree Discussion**

Now, I will show you your commentary and the discursive map made by the author of a construction at the point in the semester that they approached a task, then I will show you their proof construction and your commentary. Finally, I will show you the tree made by this student. This is the students' best attempt to describe how they thought about and completed this task. The same student is responsible for all three artifacts, and these are students that were in your classroom this semester.

1. Here is *triad discursive map* and here is *triad proof construction*. Please take a moment to compare your notes on both artifacts in this sample. Tell me any comments you have while comparing the map to the construction, as well as the rating you assigned to these artifacts.
    a. Probe: Is there anything you find surprising between the discursive map that resulted in this construction?
    b. Probe: If there anything that you find unsurprising between these two artifacts?
    c. Probe: Do you have any other comments about these artifacts and your commentary on them?

2. Here is this students' *triad tree* for this task. Recall that I am providing you with the discursive map(s) that the student stated was most representative of their function and proof discourse at the time they made this construction and tree. Tell me any comments you have while comparing the tree to the maps and construction this student created.
    a. Probe: Is there anything you find surprising between all three artifacts in this student's triad?
    b. Probe: If there anything that find unsurprising between these three artifacts?
    c. Probe: Do you have any other comments?

*Repeat for sub-protocol multiple triads.*

*Final Question:* You have examined a few constructions, trees, and discursive maps over the semester. Is there anything regarding student thinking that you are currently wondering about? Please share any commentary that you have regarding the progression of students' proof and function discourses.

*Wrap-Up:* Thank you for your time and input over the course of the semester!